\documentclass[a4paper,10pt]{article}

\usepackage[english]{babel}

\usepackage[a4paper,left=6cm,right=3cm,top=2.5cm,bottom=2.5cm]{geometry}

\usepackage{amsfonts,amsmath}
\usepackage{bm}
\usepackage{authblk}
\usepackage{siunitx}
\usepackage{xcolor}
\usepackage{subcaption}
\usepackage{float}
\usepackage{multirow}
\usepackage{comment}
\usepackage{overpic}
\usepackage{caption}
\newcommand{\fig}[1]{Fig.~\ref{#1}}{\color{blue}}

\usepackage{tikz}
\usetikzlibrary{positioning, arrows.meta, calc, shadows, shapes.geometric, shapes.symbols, fit, backgrounds, decorations.pathreplacing}

\usepackage{stackengine,graphicx}
\usepackage{tikz}
\usepackage{tabularx}
\usepackage[linesnumbered,ruled,vlined]{algorithm2e}
\usepackage[draft,inline,marginclue]{fixme}
\FXRegisterAuthor{gs}{ags}{\color{red}GS}

\usepackage{hyperref}

\usepackage{fancyhdr}
\title{Model Order Reduction of Cerebrovascular Hemodynamics Using POD–Galerkin and Reservoir Computing-based Approach}

\author[1]{Rahul Halder\footnote{rhalder@sissa.it}}

\author[1]{Arash Hajisharifi\footnote{ahajisha@sissa.it}}

\author[2,3]{Kabir Bakhshaei\footnote{ kabir.bakhshaei@santannapisa.i}}

\author[1]{Gianluigi Rozza\footnote{grozza@sissa.it}}

\affil[1]{Mathematics Area, mathLab, SISSA, via Bonomea 265, I-34136 Trieste, Italy}

\affil[2]{Biorobotics Institute, SMART Lab, Sant'Anna School of Advanced Studies, Pisa, Italy}
\affil[3]{Department of Computer Science, University of Pisa, Pisa, Italy}

\date{}

\begin{document}
\maketitle


\begin{abstract}
We investigate model order reduction (MOR) strategies for simulating unsteady hemodynamics within cerebrovascular systems, contrasting a physics-based intrusive approach with a data-driven non-intrusive framework. High-fidelity 3D Computational Fluid Dynamics (CFD) snapshots of an idealised basilar artery bifurcation are first compressed into a low-dimensional latent space using Proper Orthogonal Decomposition (POD). We evaluate the performance of a POD–Galerkin (POD–G) model, which projects the Navier–Stokes equations onto the reduced basis, against a POD–Reservoir Computing (POD–RC) model that learns the temporal evolution of coefficients through a recurrent architecture. In the present work, the introduction of a multi-harmonic and multi-amplitude training signal eliminates the requirement of multiple training signals with different amplitudes and frequencies, and thereby expedites the training process. Both methodologies achieve computational speed-ups on the order of $10^2$ to $10^3$ compared to full-order simulations, highlighting their potential as robust, efficient surrogates for real-time and
 accurate prediction of crucial flow properties such as wall shear stress (WSS), even adjacent to the bifurcation zones.

\end{abstract}


\section{Introduction}
Cerebrovascular flow simulations are frequently used in the study of intracranial aneurysms and ischemic stroke, as local hemodynamic quantities are closely linked to disease progression and rupture mechanisms \cite{Steinman2003, Cebral2005, Morris2016}. In parallel,
such simulations are also widely employed to study flow distribution and hemodynamic balance in cerebrovascular networks \cite{chnafa2017improved}. 

In practice, three-dimensional Computational Fluid Dynamics (CFD) models have been extensively used to resolve quantities such as wall shear stress (WSS) and complex flow structures that cannot be accessed through simplified analytical descriptions \cite{chnafa2017improved}.

In many applications, resolving fine-scale flow structures and multiphase interactions requires fully three-dimensional CFD simulations, as these features cannot be reliably captured by reduced or averaged models. This is the case, for example, in interface-resolved DNS of turbulent multiphase flows where the analysis relies on detailed information about flow topology and interfacial dynamics \cite{hajisharifi2021particle, hajisharifi2022interface}. Despite the physical insight they provide, the computational cost of such simulations makes them impractical for extensive parametric studies or time-sensitive applications.

Simplified lumped-parameter (0D) and one-dimensional network models are able to reproduce global flow rates and flow splits at very low computational cost \cite{chnafa2017improved}. However, by construction, these models do not provide spatially resolved information and are therefore not suitable for the analysis of local hemodynamic quantities. Their limitation is that they do not resolve the three-dimensional velocity and pressure fields, and therefore cannot be used to compute local indicators (e.g., WSS) with the required fidelity. Reduced-order models are typically introduced at this stage to reduce the cost of 3D CFD while preserving the main spatio-temporal features of the solution \cite{chnafa2017improved, macraild2024reduced, Quarteroni2016, Ballarin2016}.

A classical route to reduced modeling in fluid dynamics relies on projection-based techniques. Proper Orthogonal Decomposition (POD) is widely used to extract a small number of spatial modes that capture most of the energy contained in high-fidelity flow data \cite{benner2015survey, liu2023proper, Lumley1967, Berkooz1993, taira2017modal}. For flow regimes dominated by strong nonlinear effects or characterized by a slow Kolmogorov $n$-width decay, linear subspace approximations based on POD may provide limited compression efficiency. As a result, nonlinear manifold-based reduced-order models, often constructed using deep learning architectures such as convolutional autoencoders, have been proposed as an alternative \cite{Romor2023}.

Projecting the Navier--Stokes equations onto the reduced basis leads to a POD--Galerkin (POD--G) system that can capture the essential unsteady dynamics of the flow using a relatively small number of degrees of freedom \cite{taira2017modal, Noack2003, Stabile2018}. Within this context, stabilized projection-based reduced-order models for incompressible flows have been widely studied in finite-volume settings, with particular attention devoted to pressure stabilization techniques and the enforcement of inf--sup consistency in order to obtain reliable reduced dynamics in unsteady regimes \cite{Stabile2018}. A notable advantage of these intrusive approaches is their direct connection to the governing equations, which allows the underlying physical structure of the problem to be preserved. Nevertheless, purely projection-based ROMs may encounter difficulties when applied beyond their calibration range, particularly in the presence of strong nonlinear effects. Under these conditions, standard Galerkin formulations may experience solution drift or a gradual accumulation of truncation errors, a behavior that has prompted the investigation of more robust reduced-order modeling approaches \cite{macraild2024reduced}. In addition, the intrusive character of these methods necessitates direct access to the fully discretized equations and solver internals, which can pose practical difficulties when working with commercial solvers or large-scale, tightly coupled simulation workflows \cite{benner2015survey}.
A different strand of research moves away from equation projection altogether and instead focuses on data-driven modeling. In recent years, machine learning (ML) has been widely explored as a non-intrusive surrogate for fluid dynamics problems, especially in settings where repeated high-fidelity simulations become computationally expensive \cite{duraisamy2019turbulence, sattari2025machine, Brunton2020, Maulik2021, Pawar2020}. 
Non-intrusive POD-based reduced-order models are commonly constructed by coupling the POD representation with data-driven regression or recurrent architectures for the evolution of the reduced coefficients. Such approaches have been applied to particle-laden flows \cite{hajisharifi2023non, hajisharifi2024lstm}, conjugate heat transfer problems \cite{hajisharifi2025deep}, and mesoscale atmospheric dynamics \cite{hajisharifi2024comparison}. Since no explicit reduced governing equations are derived, these methods can be implemented without direct access to the full solver, which has contributed to their use in a range of complex flow applications.

At the same time, hybrid reduced-order modeling strategies have been proposed that sit between purely data-driven methods and classical projection-based ROMs. In these frameworks, physically motivated reduced formulations are retained, while data-driven closure models are introduced to represent unresolved energy transfer mechanisms. This combination has been shown to improve predictive accuracy compared with either purely data-driven or purely projection-based approaches \cite{IVAGNES2023127920}. Hybrid reduced-order modeling strategies have been proposed that sit between purely data-driven methods and classical projection-based ROMs, where intrusive formulations are complemented by data-driven closures to model unresolved dynamics and improve robustness and accuracy \cite{rooholamin2026hybriddiscretizethenprojectreducedorder}. Within this context, ML-augmented ROMs have received particular attention. Here, ML components are used to capture neglected dynamics or to improve numerical stability. As an example, MacRaild et al. \cite{macraild2024reduced} combined POD with neural networks for aneurysm flow simulations and reported substantial computational speed-ups.


Beyond standard forward simulations, stochastic and data-driven approaches have also found use in inverse problems, for instance in the estimation of parameters or boundary conditions in cardiovascular flow models \cite{bakhshaei2025stochastic}. In this setting, the temporal evolution of the flow is inferred directly from pre-computed solution data. Reservoir Computing (RC), a particular class of recurrent neural network, has shown promising results for chaotic and strongly unsteady systems \cite{pathak2018model, Lukosevicius2009, Pandey2020, hajisharifi2025combining}. In a different application domain, a defining feature of RC is that the nonlinear reservoir is kept fixed, and only a linear readout layer is trained. This design leads to low training cost and favorable computational efficiency, while still allowing complex temporal dependencies to be captured \cite{pathak2018model}. A related non-intrusive POD–Reservoir Computing framework has been successfully applied to industrial continuous casting processes \cite{gowrachari2025reservoir}.

\begin{figure}[H]
\centering
\begin{subfigure}[b]{1\textwidth}
\centering
{\label{fig:hemodynamics}\includegraphics[width=0.8\linewidth]{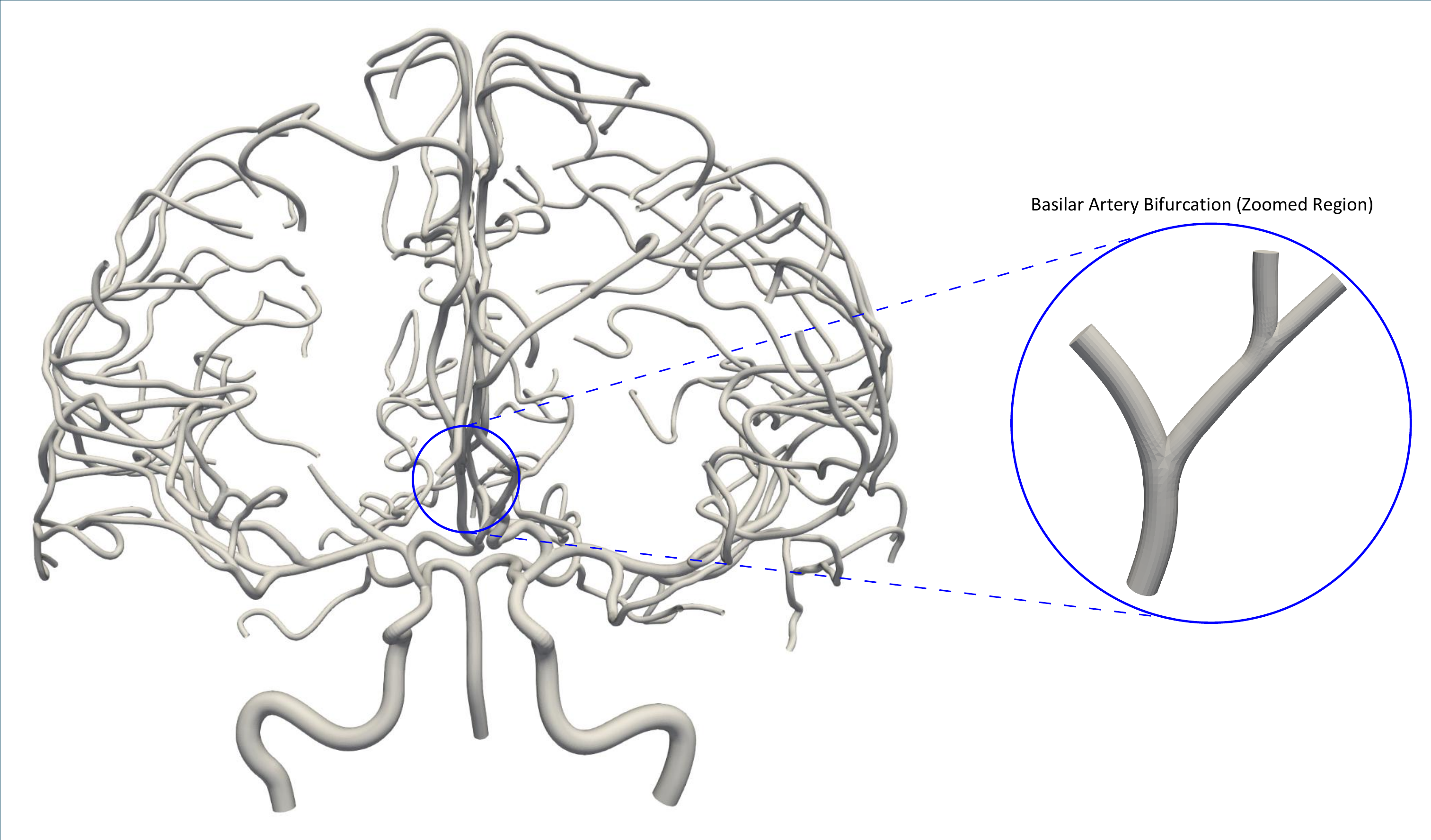}}
\end{subfigure}
\caption{
Full 3D cerebrovascular network and zoomed view of the basilar artery bifurcation.
The highlighted bifurcation region corresponds to the anatomical segment used to construct the idealized 3D computational model employed in this study.
Both the basilar artery mesh and the cerebrovascular network geometry were generated using the meshing framework of Decroocq et al.~\cite{decroocq2023meshing} and subsequently merged
for visualization purposes.}

\label{fig:ANN-DisPINN_Full}
\end{figure}

In this work, we focus on a direct comparison between these two modeling paradigms, intrusive physics-based and non-intrusive data-driven, in the context of cerebrovascular hemodynamics. An idealized basilar artery bifurcation is considered as a test case, representing the anatomical segment highlighted within the full 3D cerebrovascular network illustrated in Fig.~\ref{fig:ANN-DisPINN_Full}. We contrast an intrusive POD--Galerkin (POD--G) reduced-order model with a non-intrusive POD--Reservoir Computing (POD--RC) surrogate, both constructed from the same high-fidelity data. The comparison is carried out in terms of pressure, velocity, and WSS reconstruction, with particular attention paid to generalization across different pulsatile inflow conditions, specifically testing signals distinct from those used in training. Beyond accuracy, we assess the computational gains offered by each approach. The results indicate that both reduced models achieve speed-ups on the order of $10^2$–$10^3$ relative to full-order CFD, supporting their use as efficient surrogates for hemodynamic analysis. The contribution of the present work is highlighted as follows:

\begin{itemize}

    \item{Multiharmonic and Multiamplitude training signal employed to efficiently train both the POD-Galerkin (POD-G) and POD-Reservoir computing architecture (POD-RC).}
    \item {Reservoir Computing-based architecture is combined with Proper Orthogonal Decomposition in the context of biomedical applications and compared with the benchmark CFD results, POD-Galerkin reduced order model.} 
    \item {A comparative study is carried out between a physics-based approach and a data-driven approach for model order reduction in the context of cerebrovascular hemodynamics, and the corresponding computational time and accuracy of different approaches are thoroughly investigated. }    
\end{itemize}

\section{Governing Equations}
\label{sec:gov_equations}

The blood flow inside the cerebrovascular domain $\Omega$ is modeled as an unsteady, incompressible, laminar Newtonian fluid over the time interval $(0,\,T]$. The governing equations are the Navier--Stokes equations expressing conservation of momentum and mass:
\begin{equation}
\label{eq:NS_momentum}
\rho \frac{\partial \mathbf{u}}{\partial t} 
+ \rho (\mathbf{u} \cdot \nabla) \mathbf{u}
- \nabla \cdot \bm{\sigma} = \mathbf{0} 
\quad \text{in } \Omega \times (0, T],
\end{equation}

\begin{equation}
\label{eq:NS_continuity}
\nabla \cdot \mathbf{u} = 0 
\quad \text{in } \Omega \times (0, T],
\end{equation}

where $\rho$ denotes the fluid density, $\mathbf{u}$ is the velocity vector field, and $\bm{\sigma}$ is the Cauchy stress tensor. For laminar Newtonian flow, the stress tensor is given by:

\begin{equation}
\label{eq:stress_tensor}
\bm{\sigma} = -p \mathbf{I} + 2\mu \mathbf{D},
\end{equation}

where $p$ is the pressure, $\mathbf{I}$ is the identity tensor, $\mu$ is the constant dynamic viscosity, and $\mathbf{D}$ is the symmetric strain-rate tensor:
\begin{equation}
\label{eq:strain_rate}
\mathbf{D} = \frac{1}{2} \left( \nabla \mathbf{u} + (\nabla \mathbf{u})^T \right).
\end{equation}
Equations \eqref{eq:NS_momentum}--\eqref{eq:strain_rate} constitute the standard incompressible Navier--Stokes formulation and serve as the high-fidelity model for generating the snapshots used in the subsequent POD-based reduced-order modelling strategies.

\subsection{Laminar Flow Assumption}


Cerebral arteries operate in a low-to-moderate Reynolds-number regime because of their small diameters and physiological flow speeds \cite{chnafa2017improved}. In this setting, both intracranial and cervical flows are generally reported as laminar across the cardiac cycle, with turbulence intensity close to negligible under normal conditions \cite{chnafa2017improved, macraild2024reduced}. This is different from large-vessel flow (e.g., the aorta), where transitional features can occur around peak systole. For the cerebrovascular cases considered here, viscous and wall-shear effects dominate, and no turbulence closure is required. We therefore model the flow using the incompressible Navier--Stokes equations without a turbulence model.

\subsection{Boundary Conditions}

To simulate pulsatile cerebrovascular flow, boundary conditions are prescribed on the computational domain $\partial\Omega = \Gamma_{\mathrm{in}} \cup \Gamma_{\mathrm{out}} \cup 
\Gamma_{\mathrm{wall}}$.

\subsubsection*{Inlet boundary condition}
A time-dependent velocity inlet condition is applied at $\Gamma_{\mathrm{in}}$. 
The inlet velocity magnitude is prescribed as a multi-harmonic waveform of the form

\begin{equation}
u_{\mathrm{in}}(t)
= \bar{u}
+ \sum_{k=1}^{H}
A_k \sin\!\left( 2\pi f_k\, t + \phi_k \right),
\label{eq:inletBC}
\end{equation}

where $\bar{u} = 0.2$ is the mean inflow magnitude and the number of harmonics satisfies $H \in \{2,3,4,5\}$.  
The amplitudes, frequencies, and phases are sampled uniformly as

\[
A_k \sim \mathcal{U}(0.02,\,0.05), \qquad
f_k \sim \mathcal{U}(0.20,\,0.50), \qquad
\phi_k \sim \mathcal{U}(0,\,2\pi).
\]

This synthetic multi-harmonic waveform generates a smooth pulsatile inflow and drives the temporal dynamics of the full-order simulations,
as well as the resulting POD modal coefficients used for reduced-order modeling.

The velocity vector at the inlet is aligned with the prescribed normal direction $\mathbf{n}_{\mathrm{in}}$:

\[
\mathbf{u}(\mathbf{x},t)
= u_{\mathrm{in}}(t)\, \mathbf{n}_{\mathrm{in}}
\quad \text{on } \Gamma_{\mathrm{in}}.
\]

\begin{figure}[htbp]
    \centering
    \begin{subfigure}[b]{1\textwidth}
        \centering
        \includegraphics[width=0.8\linewidth]{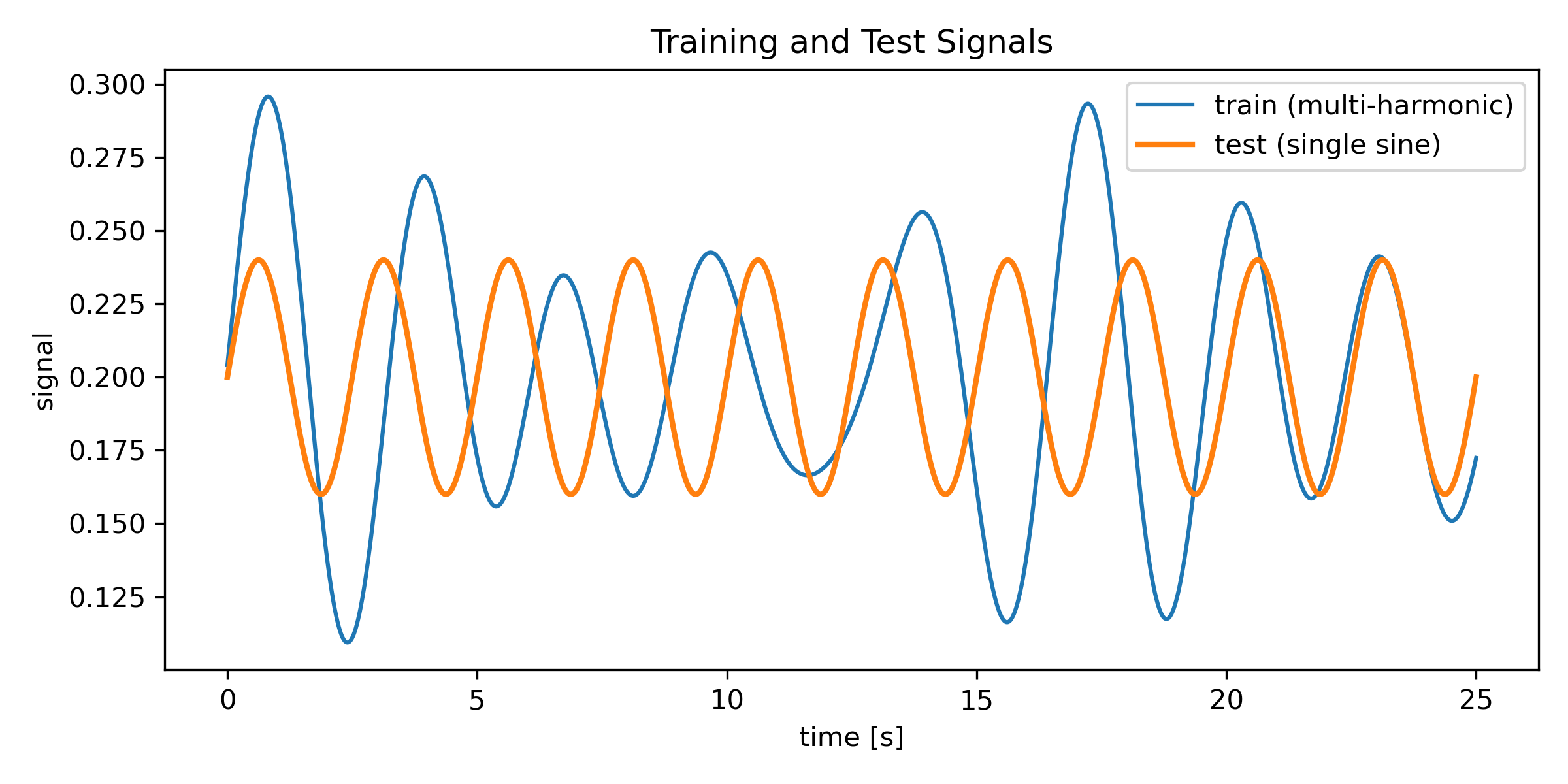}
    \end{subfigure}
    \caption{Synthetic inlet velocity boundary conditions used in this work.
The training waveform is a multi-harmonic signal, while the test waveform is a single-frequency sinusoid.
These signals drive the unsteady inflow in the FOM and are used to evaluate the temporal learning and generalization capability of the surrogate models.}
    \label{fig:training_test_signals}
\end{figure}

Figure ~\ref{fig:training_test_signals} presents the inlet velocity signals adopted in this work. The reduced-order models are trained using a smooth multi-harmonic inflow waveform, whereas their performance is evaluated using a different, single-frequency sinusoidal signal. By employing two distinct inflow conditions, we specifically test whether the proposed models can maintain accuracy when exposed to boundary conditions that were not encountered during training.

\subsubsection*{Outlet boundary condition}

At the outlet boundary $\Gamma_{\mathrm{out}}$, a traction-free (zero normal stress) condition is applied:
\begin{equation}
\bm{\sigma}(\mathbf{u},p)\, \mathbf{n} = \mathbf{0} \quad \text{on } \Gamma_{\mathrm{out}},
\end{equation}
where $\mathbf{n}$ is the outward unit normal vector. This boundary condition allows the flow to exit naturally without artificial reflections.

\subsubsection*{Wall boundary condition}

Rigid arterial walls are modeled with a no-slip condition:
\begin{equation}
\mathbf{u} = \mathbf{0} \quad \text{on } \Gamma_{\mathrm{wall}}.
\end{equation}
Together, these boundary conditions yield a well-posed problem that captures realistic pulsatile flow behavior within the cerebrovascular geometry.

For a proof of concept of ROM development, we adopt a simplified laminar Newtonian model with stress-free outlets in the present work. However, turbulence, non-Newtonian rheology, and time-varying pressure outlets (Bakhshaei et al.\ \cite{bakhshaei2025stochastic}) should be considered as future extensions of the present work.

\section{Methodology}

\begin{figure}[H]
\centering
\begin{tikzpicture}[
    node distance=1.0cm and 0.5cm,
    >=Stealth,
    font=\footnotesize\sffamily,
    process/.style={
        draw=gray!60, thick, rounded corners=2pt, fill=gray!5,
        minimum width=5.5cm, minimum height=1.0cm, align=center,
        drop shadow={opacity=0.1}
    },
    database/.style={
        cylinder, cylinder uses custom fill, shape border rotate=90,
        aspect=0.25, draw=gray!60, thick, fill=gray!10,
        minimum width=4cm, minimum height=1.2cm, align=center,
        drop shadow={opacity=0.1}
    },
    model/.style={
        draw=gray!60, thick, rounded corners=2pt,
        minimum width=4.0cm, minimum height=1.2cm, align=center,
        drop shadow={opacity=0.15}
    },
    connector/.style={
        ->, thick, color=gray!70, rounded corners=3pt
    }
]


    \node[process] (fom) {
        \textbf{High-Fidelity CFD (FOM)}\\
        Full-order Navier--Stokes simulation\\
        3D basilar artery bifurcation
    };

    \node[database, below=0.8cm of fom] (snap) {
        \textbf{Snapshot Collection}\\
        Velocity \& pressure fields\\
        $\mathbf{S} = [\mathbf{u}_1, \dots, \mathbf{u}_{N_s}]$
    };

    \node[process, below=0.8cm of snap] (pod) {
        \textbf{POD Decomposition}\\
        Modes $\Phi(\mathbf{x})$, coefficients $a(t)$\\
        Surrogates predict $a(t)$
    };

    \coordinate (modelsCenter) at ($(pod.south)+(0,-1.7cm)$);

    \node[model, fill=blue!15] (rom) at ($(modelsCenter)+(-2.6cm,0)$) {
        \textbf{POD--Galerkin ROM}\\
        (Intrusive, physics-based)
    };

    \node[model, fill=orange!15] (rc) at ($(modelsCenter)+( 2.6cm,0)$) {
        \textbf{Reservoir Computing (RC)}\\
        (Non-intrusive)\\
        trained on $a(t)$
    };

    \node[process, minimum width=8cm] (compare) at ($(modelsCenter)+(0,-2.2cm)$) {
        \textbf{Field Reconstruction \& Model Evaluation}\\
        Reconstruct $\tilde{\mathbf{u}}(\mathbf{x},t) \approx \sum_{i=1}^{k} \tilde{a}_i(t)\,\Phi_i(\mathbf{x})$ and compare with FOM\\
        Compute error metrics w.r.t.\ FOM data
    };


    \draw[connector] (fom) -- (snap);
    \draw[connector] (snap) -- (pod);

    \draw[thick, color=gray!70] (pod.south) -- ++(0,-0.6) coordinate (split);
    \draw[connector] (split) -| (rom.north);
    \draw[connector] (split) -| (rc.north);

    \path (rom.south) -- ++(0,-0.6) coordinate (romDown);
    \path (rc.south)  -- ++(0,-0.6) coordinate (rcDown);
    \path (romDown) -- (rcDown) coordinate[midway] (merge);

    \draw[thick, color=gray!70] (rom.south) -- (romDown);
    \draw[thick, color=gray!70] (rc.south)  -- (rcDown);
    \draw[thick, color=gray!70] (romDown) -- (merge);
    \draw[thick, color=gray!70] (rcDown)  -- (merge);
    \draw[connector] (merge) -- (compare.north);

    \begin{pgfonlayer}{background}
        \node[
            fit=(fom)(snap)(pod),
            fill=gray!2, draw=gray!30, dashed, rounded corners=5pt,
            inner sep=10pt, yshift=2pt
        ] (offline) {};
        \node[anchor=north east, text=gray!80, font=\bfseries\scriptsize]
            at (offline.north east) {OFFLINE PHASE};

        \node[
            fit=(rom)(rc)(compare),
            fill=blue!2, draw=blue!10, dashed, rounded corners=5pt,
            inner sep=10pt, yshift=-2pt
        ] (online) {};
        \node[anchor=south east, text=blue!40, font=\bfseries\scriptsize]
            at (online.south east) {ONLINE PHASE};
    \end{pgfonlayer}

\end{tikzpicture}
\caption{Schematic of the surrogate workflow. The Offline Phase (gray) consists of FOM data generation and POD basis extraction. In the Online Phase (blue), the temporal coefficients are predicted either by an intrusive POD--Galerkin ROM or by a non-intrusive POD--RC surrogate; reconstructed fields are then evaluated against the FOM reference.}
\label{fig:flowchart_improved}
\end{figure}
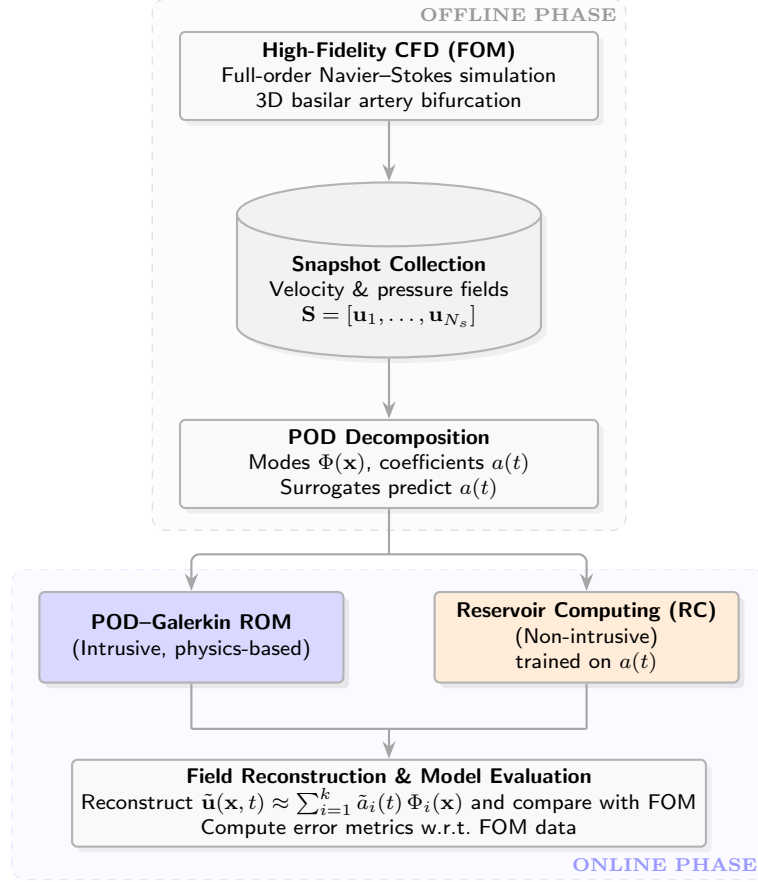

In this section, we will briefly introduce two different model order reduction strategies, one of which is data-driven, and the other one is a physics-driven and projection-based approach to expedite the numerical simulation associated with cerebrovascular hemodynamics to ensure real-time and reliable prediction of flow physics for patient-specific diagnosis.  Figure~\ref{fig:flowchart_improved} summarizes the workflow adopted in this work. We separate the operations that are computationally demanding from those that must run quickly at inference time. For clarity, the figure groups these steps into an Offline Phase (gray) and an Online Phase (blue).

\vspace{0.3cm}

\noindent \textbf{Offline Phase:} We start by running high-fidelity CFD simulations (FOM) on the 3D basilar artery bifurcation. The resulting velocity and pressure solutions corresponding to the training signal in \ref{fig:training_test_signals} are first collected into a snapshot matrix $\mathbf{S}$. Proper Orthogonal Decomposition (POD) is then performed on $\mathbf{S}$ to obtain a reduced basis of spatial modes $\Phi(\mathbf{x})$ together with the corresponding time coefficients $a(t)$. After this step, the full 3D fields are represented through a small set of modal amplitudes.

\vspace{0.3cm}

\noindent \textbf{Online Phase:} The goal is to advance the coefficients in time, i.e., to predict ${a}(t)$. We test two alternatives:

\begin{enumerate}
    \item \textbf{POD--Galerkin ROM (Intrusive):} we derive a reduced dynamical system by projecting the incompressible Navier--Stokes equations into the POD subspace.
    \item \textbf{POD--Reservoir Computing (RC) ROM (Non-intrusive):} we learn the time evolution directly from the coefficient data using a reservoir computing model.
\end{enumerate}


In both cases, the surrogate produces $\tilde{a}(t)$, an approximation of the temporal coefficient $a(t)$ in physics-based and data-driven ways. The corresponding 3D velocity and pressure fields are reconstructed via the POD expansion and are subsequently compared with the FOM solutions to compute the error metrics reported in this study. Therefore, we will first demonstrate the Proper Orthogonal Decomposition (POD) approach, which is the common step for both methodologies.

\subsection{Snapshot Collection and Proper Orthogonal Decomposition}
\label{sec:intrusive_pod}

Here, we first discuss snapshot collection of different solution fields and application of proper orthogonal decomposition, which are common steps for both the physics-based and data-driven approaches. To carry out the steps, we first assume that the full order model can be expressed as a linear combination of the spatial modes and their corresponding temporal coefficients are then computed subsequently as mentioned by \cite{Stabile2018, star2019novel}. The velocity snapshots $\boldsymbol{u}\left(\boldsymbol{x}, t_n\right)$ and pressure snapshots $p\left(\boldsymbol{x}, t_n\right)$ at time $t_n$ are computed, respectively, by the following relationships, where $\varphi_i$ and $\chi_i$ are the modes for velocity and pressure, respectively.  

\begin{equation}\label{eq:grad_burgers_main}
\begin{aligned}
\boldsymbol{u}\left(\boldsymbol{x}, t_n\right) \approx \boldsymbol{u}_r\left(\boldsymbol{x}, t_n\right)=\sum_{i=1}^{N_r^u} \boldsymbol{\varphi}_i(\boldsymbol{x}) a_i\left(t_n\right), \\
\quad p\left(\boldsymbol{x}, t_n\right) \approx p_r\left(\boldsymbol{x}, t_n\right)=\sum_{i=1}^{N_r^p} \chi_i(\boldsymbol{x}) b_i\left(t_n\right),
\end{aligned}
\end{equation}

$\boldsymbol{a}\left(t_n\right)=\left[a_1\left(t_n\right),\ldots, a_n\left(t_n\right)\right]^T$ and $\boldsymbol{b}\left(t_n\right)=\left[b_1\left(t_n\right), \ldots, b_n\left(t_n\right)\right]^T$ are column vectors which contain the corresponding time-dependent coefficients. $N_r^u$ is the number of velocity modes, and $N_r^p$ is the number of pressure modes. Furthermore, the modes are orthonormal to each other: $\left(\varphi_i, \varphi_j\right)_{L^2(\Omega)}=\delta_{i j}$, where $\delta$ is the Kronecker delta.  We will now briefly introduce the Proper Orthogonal Decomposition (POD) approach, followed by the data-driven approach using the POD-RC network and the physics-based approach, POD-Galerkin, which utilises the projection of the incompressible Navier-Stokes equation onto the reduced basis. The optimal POD basis space for velocity, $E_u^{POD}=\left[\varphi_1, \varphi_2, \ldots, \varphi_{N_r^u}\right]$, is then constructed by minimizing the difference between the snapshots and their orthogonal projection onto the basis for the $L^2$-norm \cite{Lumley1967,Sirovich1987,Kunisch2002,Benner2015, Stabile2018, star2019novel}. This gives the following minimization problem. First, we demonstrate the basics of the proper orthogonal decomposition approach considered in the present work:

\begin{equation}\label{eq:grad_burgers_main}
\begin{aligned}
E_u^{P O D} =\arg \min _{\boldsymbol{\varphi}_1, \ldots, \boldsymbol{\varphi}_{N_s^u}} \frac{1}{N_s^u} \sum_{n=1}^{N_r^u}\left\|\boldsymbol{u}\left(\boldsymbol{x}, t_n\right)-\sum_{i=1}^{N_r^u}\left(\boldsymbol{u}\left(\boldsymbol{x}, t_n\right), \boldsymbol{\varphi}_{\boldsymbol{i}}(\boldsymbol{x})\right)_{L^2(\Omega)} \boldsymbol{\varphi}_i(\boldsymbol{x})\right\|_{L^2(\Omega)}^2,
\end{aligned}
\end{equation}

where, $N_s^u$ is the number of collected velocity snapshots from training and $N_r^u$ is the dimension of the POD spaces $\text { with } 1 \leq N_r^u \leq N_s^u$. Now, POD modes are computed by solving the following eigenvalue problem on the snapshots: 

\begin{equation}\label{eq:grad_burgers_main}
\begin{aligned}
C Q=Q \lambda
\end{aligned}
\end{equation}

where, $C_{i j}=\left(\boldsymbol{u}\left(\boldsymbol{x}, t_i\right), \boldsymbol{u}\left(\boldsymbol{x}, t_j\right)\right)_{L^2(\Omega)}$ for $i, j=1, \ldots, N_s^u$ is the correlation matrix. $Q \in \mathbb{R}^{N_s^u \times N_s^u}$ is a square matrix of eigen vectors and $\boldsymbol{\lambda} \in \mathbb{R}^{N_s^u \times N_s^u}$ is a diagonal matrix containing the eigen values. The POD modes, $\varphi_i$ are then computed as follows:

\begin{equation}\label{eq:grad_burgers_main}
\begin{aligned}
\boldsymbol{\varphi}_i(\boldsymbol{x})=\frac{1}{N_s^u \sqrt{\lambda_i}} \sum_{n=1}^{N_s^u} \boldsymbol{u}\left(\boldsymbol{x}, t_n\right) Q_{i n} \quad \text { for } i=1, \ldots, N_r^u
\end{aligned}
\end{equation}

The procedure is the same for obtaining the pressure modes. Once the POD modes are computed governing Navier-Stokes equations are projected on the reduced basis.

\subsection{Intrusive ROM: POD-Galerkin Projection}

In this section, we describe the POD-Galerkin projection-based reduced-order model. Projection-based reduced-order models formulated within a finite-volume framework have been successfully applied to canonical unsteady incompressible flow problems. The intrusive-reduced order model can be obtained by projecting the Navier-Stokes equation on POD spaces using Galerkin and Petrov-Galerkin projection as follows, by projecting the momentum equation on velocity basis vectors ($\phi$) and the continuity equation on the pressure spaces ($\chi$). 

\begin{equation}\label{eq:projection}
\begin{aligned}
\left\langle\phi_i, \frac{\partial \boldsymbol{u}}{\partial t}+\nabla \cdot(\boldsymbol{u} \otimes \boldsymbol{u})-\nabla \cdot(\nu \nabla \boldsymbol{u})+\nabla p\right\rangle=0 \quad \text { in } \Omega \times[0, T] . \\
\left\langle\chi_i, \nabla \cdot \boldsymbol{u}\right\rangle=0, \text { in } \Omega \times[0, T] .
\end{aligned}
\end{equation}

After the projection of both equations on the reduced basis, we find the following ODEs:
\begin{equation}\label{eq:ODE-classical}
\begin{aligned}
\sum_{j=1}^{N_u^r} M_{i j} \frac{\partial \alpha_j}{\partial t} & =\sum_{j=1}^{N_u^r} \sum_{k=1}^{N_u^r} Q_{i j k} \alpha_j \alpha_k+\nu \sum_{i=1}^{N_u^r} L_{i j} \alpha_i-\sum_{i=1}^{N_p^r} P_{i j} b_i \\
\sum_{j=1}^{N_p^r} R_{i j} \alpha_j & =0
\end{aligned}
\end{equation}

where the matrices associated with the ODEs, i.e., ($M$, $Q$, $L$, $P$ and $R$) represent the reduced-order operator. They correspond to mass, nonlinear convection, diffusion, pressure gradient coupling, and incompressibility constraints. These operators can be precomputed offline, which leads to an efficient online stage. They can be written as follows: 

\begin{equation}\label{eq:Matrices}
\begin{aligned}
M_{i j} & =\left\langle\boldsymbol{\phi}_i, \boldsymbol{\phi}_j\right\rangle_{L_2(\Omega)}, \\
Q_{i j k} & =\left\langle\nabla \cdot\left(\boldsymbol{\phi}_i \otimes \boldsymbol{\phi}_j\right), \boldsymbol{\phi}_k\right\rangle_{L_2(\Omega)}, \\
L_{i j} & =\left\langle\nu \Delta \boldsymbol{\phi}_i, \boldsymbol{\phi}_j\right\rangle_{L_2(\Omega)}, \\
P_{i j} & =\left\langle\nabla \psi_i, \boldsymbol{\phi}_j\right\rangle_{L_2(\Omega)}, \\
R_{i j} & =\left\langle\nabla \cdot \boldsymbol{\phi}_i, \psi_j\right\rangle_{L_2(\Omega)} .
\end{aligned}
\end{equation}

Standard Galerkin projection-based reduced-order models are unreliable when applied to the nonlinear, unsteady Navier–Stokes equations. For fluid problems that are solved numerically using a Finite Volume discretization technique, often a Poisson Equation is solved for pressure.
As there is no dedicated equation for pressure in the Equation, the pressure Poisson equation is obtained by taking the divergence of the momentum equation and exploiting the divergence-free constraint $\nabla \cdot \boldsymbol{u}=\mathbf{0}$. The governing equation associated with the pressure Poisson equation can be written as follows:

\begin{equation}
\begin{aligned}
\frac{\partial \boldsymbol{u}}{\partial t}
+ \nabla \cdot (\boldsymbol{u} \otimes \boldsymbol{u})
- \nabla \cdot (\nu \nabla \boldsymbol{u})
= - \nabla p && \text{in } \Omega \times [0, T], \\[0.5em]
\Delta p
= - \nabla \cdot \bigl(\nabla \cdot (\boldsymbol{u} \otimes \boldsymbol{u})\bigr)
&& \text{in } \Omega \times [0, T], \\[0.5em]
\boldsymbol{n} \cdot \nabla p
= - \boldsymbol{n} \cdot \left( \nu \nabla \times \nabla \times \boldsymbol{u}
+ \frac{\partial u_{\mathrm{in}}}{\partial t} \right)
&& \text{on } \Gamma \times [0, T].
\end{aligned}
\end{equation}

Now, we intend to project this coupled system onto the reduced velocity space spanned by $\varphi_i$, and the reduced pressure space spanned by $\chi_i$, which leads to the following weak formulation:

\begin{equation}
\begin{aligned}
&\left\langle\varphi_{\boldsymbol{i}}, \boldsymbol{u}_{\boldsymbol{t}}+\boldsymbol{\nabla} \cdot(\boldsymbol{u} \otimes \boldsymbol{u})+\boldsymbol{\nabla} p-\boldsymbol{\nabla} \cdot 2 \nu \boldsymbol{\nabla}^{\boldsymbol{s}} \boldsymbol{u}\right\rangle_{L_2(\Omega)}=0,\\
& \left\langle\nabla \chi_i, \nabla p\right\rangle_{L_2(\Omega)}+\left\langle\nabla \chi_i, \nabla \cdot(\boldsymbol{u} \otimes \boldsymbol{u})\right\rangle_{L_2(\Omega)} \\
& -\nu\left\langle\boldsymbol{n} \times \nabla \chi_i, \nabla \times \boldsymbol{u}\right\rangle_{\Gamma}-\left\langle\chi_i, \boldsymbol{n} \cdot u_{\mathrm{in}}\right\rangle_{\Gamma}=0 .
\end{aligned}
\end{equation}

which further leads to reduced order formulation as follows:

\begin{equation}   
\begin{aligned}
& \boldsymbol{M}_{\boldsymbol{r}} \dot{\boldsymbol{a}}-\nu \boldsymbol{A}_{\boldsymbol{r}} \boldsymbol{a}+\boldsymbol{a}^T \mathbf{C}_{\boldsymbol{r}} \boldsymbol{a}+\boldsymbol{B}_{\boldsymbol{r}} \boldsymbol{b}=0 \\
& \boldsymbol{D}_{\boldsymbol{r}} \boldsymbol{b}+\boldsymbol{a}^T \mathbf{G}_{\boldsymbol{r}} \boldsymbol{a}-\nu \boldsymbol{N}_{\boldsymbol{r}} \boldsymbol{a}-\boldsymbol{F}_{\boldsymbol{r}}=0
\end{aligned}
\end{equation}

These additional operators introduced in this formulation can be defined as follows:

\begin{equation}    
\begin{aligned}
& D_{r_{i j}}=\left\langle\boldsymbol{\nabla} \chi_i, \boldsymbol{\nabla} \chi_j\right\rangle_{L_2(\Omega)}, \\
& \mathrm{G}_{\mathrm{r} i j k}=\left\langle\nabla \chi_i, \boldsymbol{\nabla} \cdot\left(\boldsymbol{\varphi}_{\boldsymbol{j}} \otimes \boldsymbol{\varphi}_{\boldsymbol{k}}\right)\right\rangle_{L_2(\Omega)}, \\
& N_{r_{i j}}=\left\langle\boldsymbol{n} \times \boldsymbol{\nabla} \chi_i, \boldsymbol{\nabla} \times \boldsymbol{\varphi}_{\boldsymbol{j}}\right\rangle_{\Gamma}, \\
& F_{r_i}=\left\langle\chi_i, \boldsymbol{n} \cdot u_{\mathrm{in}}\right\rangle_{\Gamma} .
\end{aligned}
\end{equation}

This formulation provides a consistent treatment of the pressure field, which leads to improved numerical robustness as compared to the classical velocity pressure POD-Galerkin formulation demonstrated in \ref{eq:ODE-classical}

\subsubsection{Boundary Conditions and Snapshot Homogenization}

In POD-Galerkin-based reduced-order models, the treatment of non-homogeneous boundary conditions is an important step, as mentioned in \cite{star2019novel}. Since POD modes are constructed as linear combinations of snapshot solutions, the resulting basis functions inherit the boundary behaviour present in the training data. Therefore, the reduced-order approximation fails to satisfy prescribed time-dependent Dirichlet boundary conditions. It can negatively affect both accuracy and stability.
Two commonly used strategies exist to address this issue: the lifting function method and the penalty approach \cite{star2019novel}. In this work, we employ the lifting function, which is simple yet based on a strong theoretical foundation. The main idea is to decompose the velocity field into a sum of a particular function that satisfies the non-homogeneous boundary conditions and a homogeneous fluctuation field. POD is then applied only to the homogeneous component. The velocity snapshots are homogenized as follows:

\begin{equation}  
\boldsymbol{u}^{\prime}(\boldsymbol{x}, t)=\boldsymbol{u}(\boldsymbol{x}, t)-\sum_{j=1}^{N_{B C}} \boldsymbol{\zeta}_{c_j}(\boldsymbol{x}) u_{B C_j}(t),
\end{equation}

where $N_{B C}$ is the number of non-homogeneous $\mathrm{BCs}, \boldsymbol{\zeta}_c(\boldsymbol{x})$ the normalized lifting functions and $u_{B C}$ is the normalized value of the corresponding Dirichlet boundary condition. The detailed approach is demonstrated in \cite{star2019novel}.

\subsection{Non-Intrusive ROM: POD-Reservoir Computing Architecture}
\label{sec:POD-RC}

Reservoir computing (RC) is a machine learning framework which efficiently captures temporal dynamics. This specific class of learning approach falls under the category of Recurrent Neural Network (RNN) from the viewpoint of machine learning. It offers an effective alternative to traditional deep-learning methods, such as Long Short Term Memory based network, as shown by Halder et al. \cite{halder2025physics,halder2020deep,halder2020deep}, Transformer based network \cite{solera2024beta}, which involve a long training time. A reservoir computing framework, on the contrary, exploits a fixed, randomly initialized dynamical system as a reservoir, and the output layer of the reservoir computing is only trained by a computationally cheap regularized linear least-squares optimization procedure, thereby drastically reducing the training time as compared to a conventional deep-learning approach. The typical RC consists of three components:  

\begin{itemize}
    \item Input Layer: This layer maps transient input to the reservoir. 

    \item Reservoir Layer: This layer is fixed and randomly initialized before training, and the input data is transformed into a complex temporal pattern.

    \item Output Layer: This is the only trainable layer in the RC, which transforms the reservoir state to the output layer, typically using linear regression techniques.
\end{itemize}

\begin{figure}[htbp]
    \centering
    \begin{subfigure}[b]{1\textwidth}
        \centering
        \includegraphics[width=0.8\linewidth]{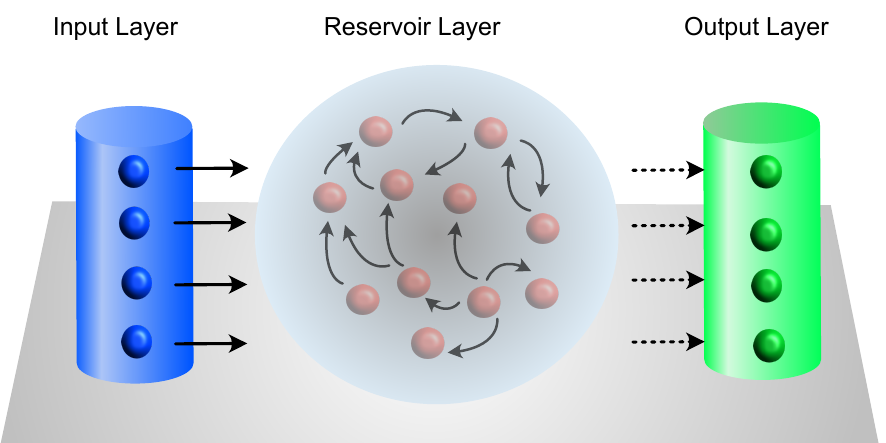}
    \end{subfigure}
    \caption{Layers of Reservoir Computing Architecture, first Layer is the input layer, second one is fixed and randomly initialized, and only the final output layer is trainable.}
    \label{fig:training_test_signals}
\end{figure}

There are three different types of reservoir computing networks. Out of which, in the present work, the Echo State Network (ESN) architecture is used. The key idea is to use a fixed, high-dimensional nonlinear dynamical system to transform inputs into a rich feature space, while training only a linear readout layer. We employ an RC network to learn the nonlinear mapping
\begin{equation}
\boldsymbol{I}:\; \bm{u}_{\text{in},k} \in \mathbb{R}^{d_{\text{in}}}
\;\longmapsto\;
\boldsymbol{\Theta}_k \in \mathbb{R}^{d_{\text{out}}},
\end{equation}
where the input $\boldsymbol{I}:\bm{u}_{\text{in},k}$ corresponds to the inlet velocity waveform of the cerebrovascular system, and the output $\boldsymbol{\Theta}_k$ represents the POD coefficients associated with pressure $p$ and velocity $v$.

Let $\bm{u}_{\text{in},k} \in \mathbb{R}^{d_{\text{in}}}$ denote the inlet velocity input at discrete time step $k$, and let $\bm{\Theta}_k \in \mathbb{R}^{d_{\text{out}}}$ denote the desired output. The reservoir state $\bm{r}_k \in \mathbb{R}^{N_r}$ evolves according to
\begin{equation}
\mathbf{r}(t+1)=\underbrace{(1-\gamma)\mathbf{r}(t)}_{\text{Linear memory}} 
+ \underbrace{\gamma f\left[\beta\left(\varepsilon {W}_{\text{in}} \bm{u}_{\text{in}}(t)
+ {W}\mathbf{r}(t) + \eta \mathbf{b}\right)\right]}_{\text{Nonlinear activation}},
\label{eq:reservoir-update}
\end{equation}
where $W_{\text{in}} \in \mathbb{R}^{N_r \times d_{\text{in}}}$ maps the inlet velocity input to the reservoir, $W \in \mathbb{R}^{N_r \times N_r}$ is the recurrent weight matrix (typically sparse), $\gamma \in (0,1]$ is the leaking rate, and $f(\cdot)$ is a nonlinear activation function such as $\tanh$.

The output is obtained through a linear readout
\begin{equation}
\bm{\Theta}_k = W_{\text{out}} \bm{r}_k + \bm{b}_{\text{out}},
\label{eq:readout}
\end{equation}
where $W_{\text{out}} \in \mathbb{R}^{d_{\text{out}} \times N_r}$ and $\bm{b}_{\text{out}} \in \mathbb{R}^{d_{\text{out}}}$ are the only trainable parameters.

During training, the reservoir states are collected in the matrix
\[
R =
\begin{bmatrix}
\bm{r}_1^\top \\
\bm{r}_2^\top \\
\vdots \\
\bm{r}_T^\top
\end{bmatrix},
\qquad
Y =
\begin{bmatrix}
\bm{\Theta}_1^\top \\
\bm{\Theta}_2^\top \\
\vdots \\
\bm{\Theta}_T^\top
\end{bmatrix},
\]
and the readout matrix is obtained via ridge regression
\begin{equation}
W_{\text{out}} =
\arg\min_W \| R W^\top - Y \|_F^2 + \lambda \| W \|_F^2,
\label{eq:ridge}
\end{equation}
where $\lambda \ge 0$ is a regularization parameter. In this section, we have demonstrated that, for a specific multi-harmonic signal as a velocity input, first, a full-order model is generated, followed by the application of Proper-Orthogonal Decomposition (POD) on the velocity and pressure snapshots. Finally, the POD coefficients are computed using a physics-based Galerkin projection and a data-driven approach, i.e., Reservoir Computing, and the associated numerical results will be discussed in the next section.


\section{Results}

In this section, we evaluate the efficiency of the POD–Galerkin reduced-order model (POD-G) and the POD–Reservoir Computing model (POD-RC) on a benchmark consisting of an idealized cerebrovascular bifurcation. This geometry has the usual hemodynamic characteristics of intracranial arterial systems, including flow division, recirculation areas, and asymmetric velocity distributions. The ROMs were constructed from high-fidelity CFD snapshots. We assess the accuracy, temporal stability and reconstruction quality of their prediction by comparing them with the reference solution.

The section is organized as follows. The cerebrovascular geometry and the resulting POD decomposition are described in Sec ~\ref{sec:POD Decomposition}. The results obtained from POD-Galerkin projection ROM and POD-Reservoir computing ROM are presented in Sec~\ref{sec:POD-G} and Sec~\ref{sec:POD-G}, respectively, including their quantitative error matrices and qualitative field reconstruction. A direct comparison between the two ROMs is provided in Sec~\ref{sec:comparison} highlighting their accuracy and robustness. Finally, the computational efficiency of the two approaches is discussed in Sec~\ref{sec:computational_cost}.


\subsection{Cerebrovascular Artery and POD Decomposition}
\label{sec:POD Decomposition}

\begin{figure}[htb]
    \centering
    \begin{subfigure}[b]{1\textwidth}
        \centering
        \subfloat[Modes = $1$]{\label{fig:umode1}\includegraphics[width=.33\linewidth]{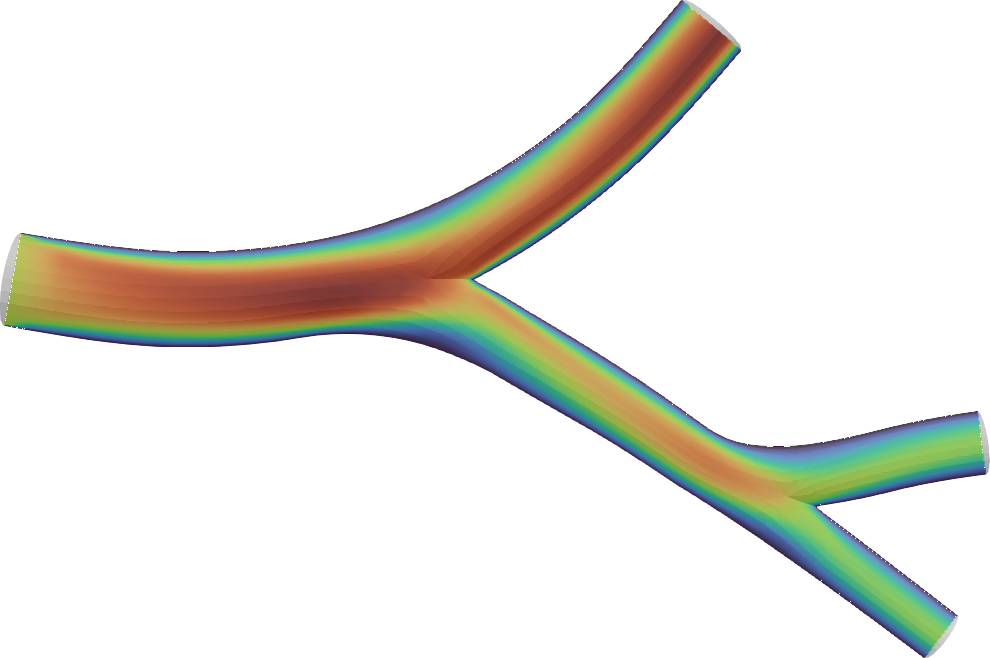}}
        \subfloat[Modes = $2$]{\label{fig:umodes2}\includegraphics[width=.33\linewidth]{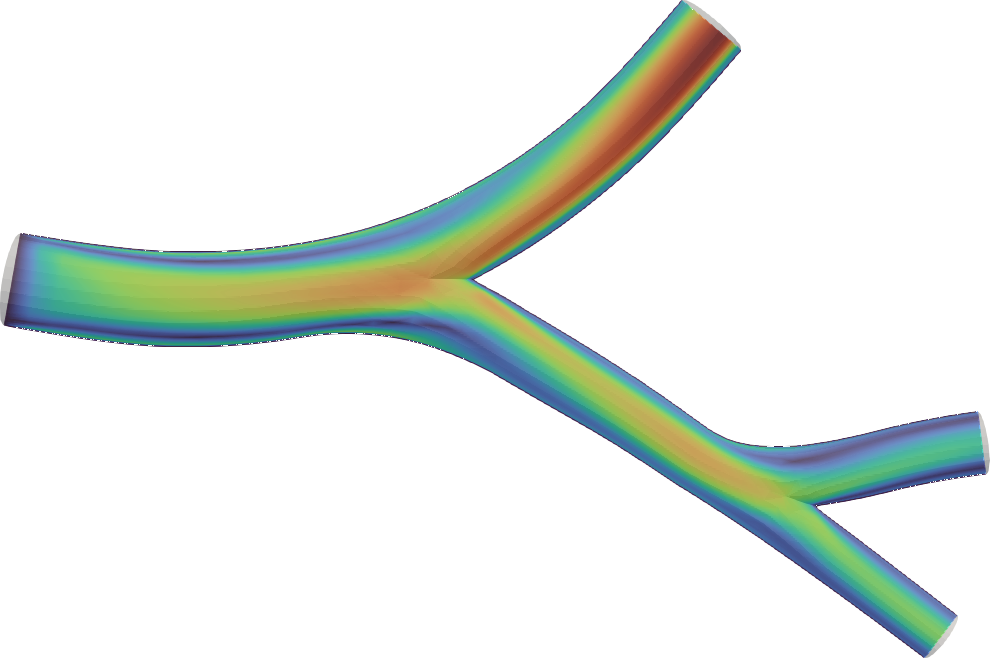}}
        \subfloat[Modes = $3$]{\label{fig:umodes3}\includegraphics[width=.33\linewidth]{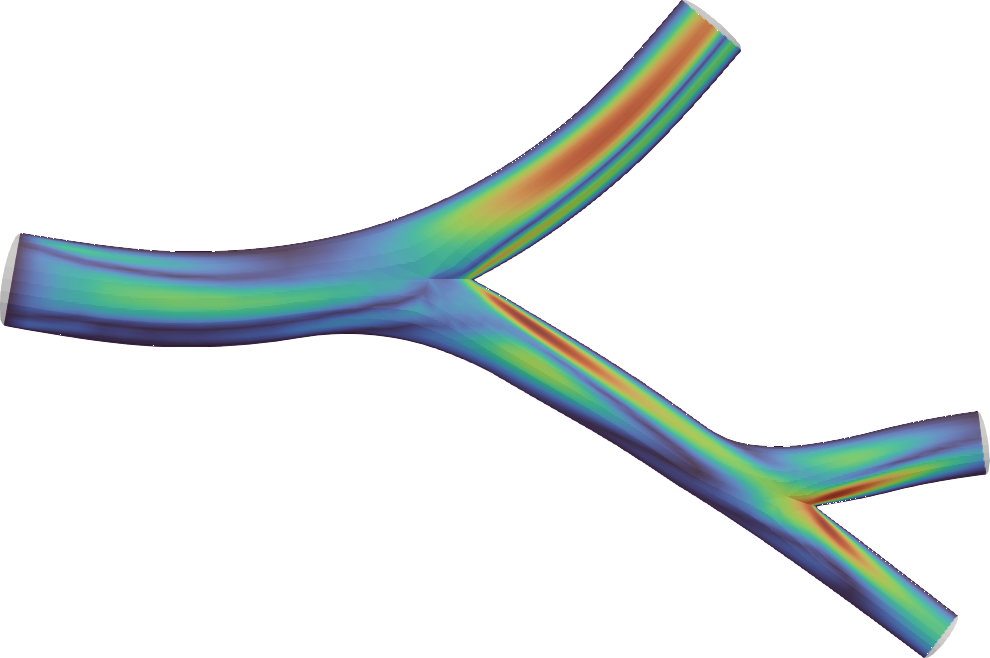}}
    \end{subfigure}

    \vspace{0.3cm} 

    \begin{subfigure}[b]{1\textwidth}
        \centering
        \subfloat[Modes = $4$]{\label{fig:umodes4}\includegraphics[width=.33\linewidth]{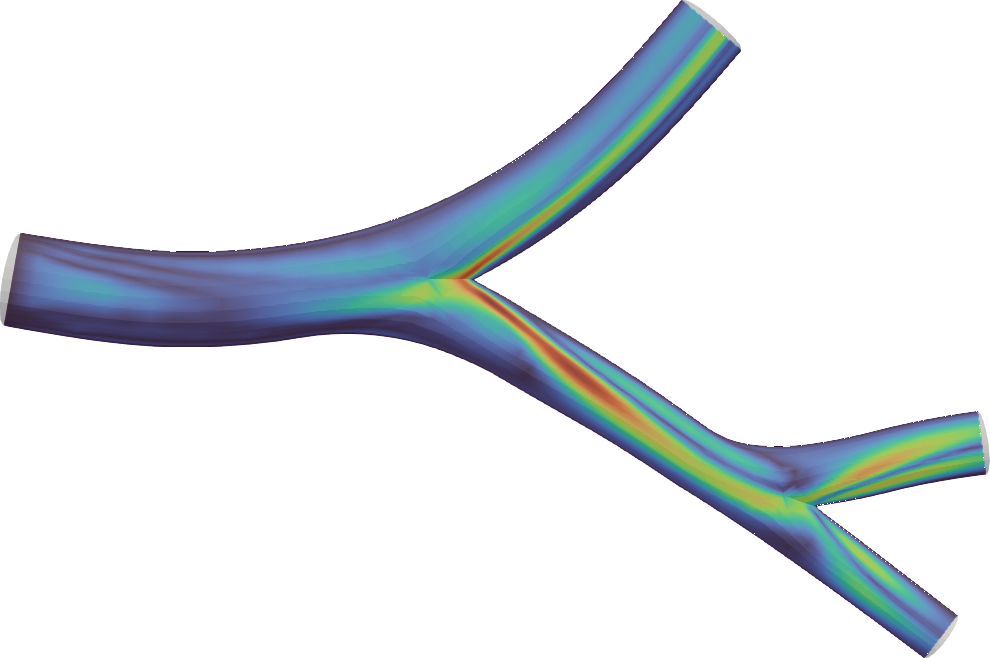}}
        \subfloat[Modes = $5$]{\label{fig:umodes5}\includegraphics[width=.33\linewidth]{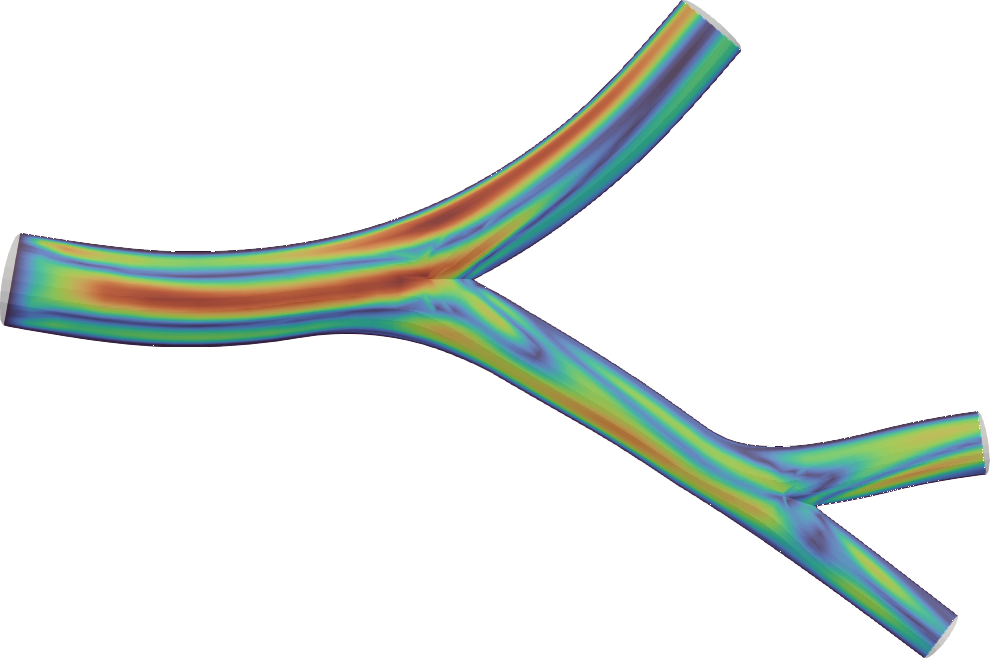}}
        \subfloat[Modes = $6$]{\label{fig:umodes6}\includegraphics[width=.33\linewidth]{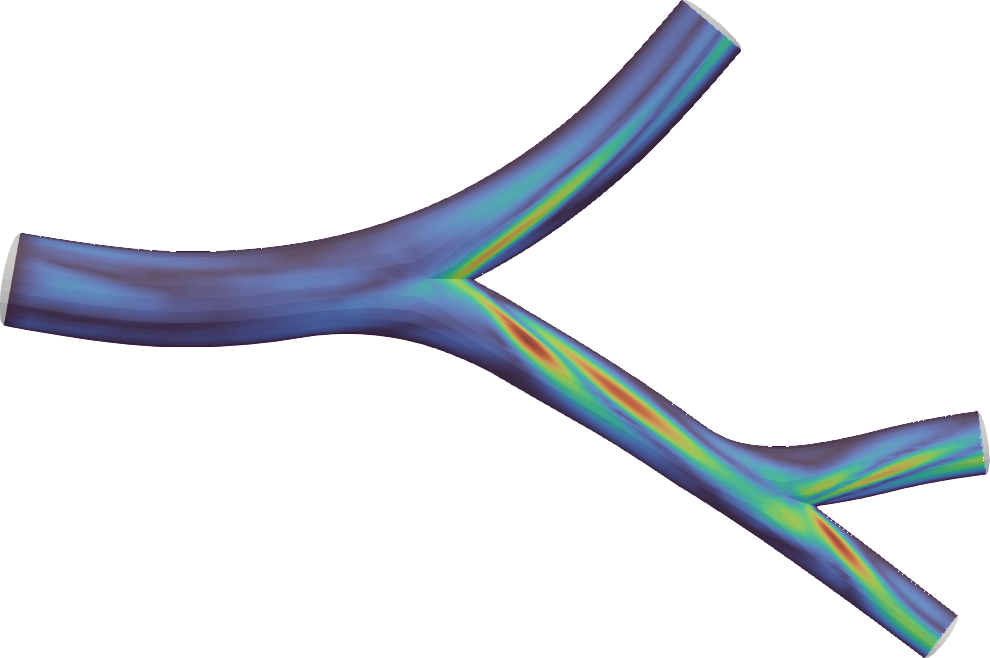}}
    \end{subfigure}

    \begin{subfigure}[b]{1\textwidth}
        \centering
        {\label{fig:uscale1}\includegraphics[width=0.25\linewidth]{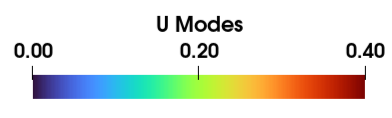}}          
    \end{subfigure}

    \medskip
    \caption{ First six POD modes of the velocity magnitude field for the cerebrovascular bifurcation geometry. The modes are ordered according to decreasing energetic contribution with the first mode representing the dominant flow structure.
}
    \label{fig:POD_U_modes}
\end{figure}

The considered benchmark in this study is an idealized cerebrovascular bifurcation, representing the basilar artery geometry. This setup represents the mainly observed hemodynamic features in intracranial arterial systems such as the flow division at the bifurcation, velocity redistribution between the daughter branches and recirculation regions. Proper Orthogonal Decomposition (POD) is employed to reduce the high-dimensional full-order CFD simulation to its reduced representation by extracting the dominant spatial structure of the flow.

Figure~\ref{fig:POD_U_modes} shows the first six POD modes of the velocity magnitude field, ordered based on their energetic contribution. The first mode represents the most dominant spatial structure, namely the mean flow topology of the bifurcation. It captures the main flow pattern in the parent vessel and its subsequent redistribution into the daughter branches.  As the number of POD modes increases, finer spatial structure and more localized variations of the velocity field are represented. In particular, the intermediate  modes shown in \fig{fig:POD_U_modes} capture asymmetries between the branches and variations near the bifurcation region, while the higher number of modes represent more complex and localized spatial structure, corresponding to the smaller-scale flow features. These results highlight the capability of POD decomposition to separate the dominant large-scale flow feature from the higher-frequency spatial structure,  compressing efficiently the unsteady flow field to its lower-dimensional representation.

\begin{figure}[t]
    \centering
    \begin{subfigure}[b]{1\textwidth}
        \centering
        \subfloat[Modes = $1$]{\label{fig:pmode1}\includegraphics[width=.33\linewidth]{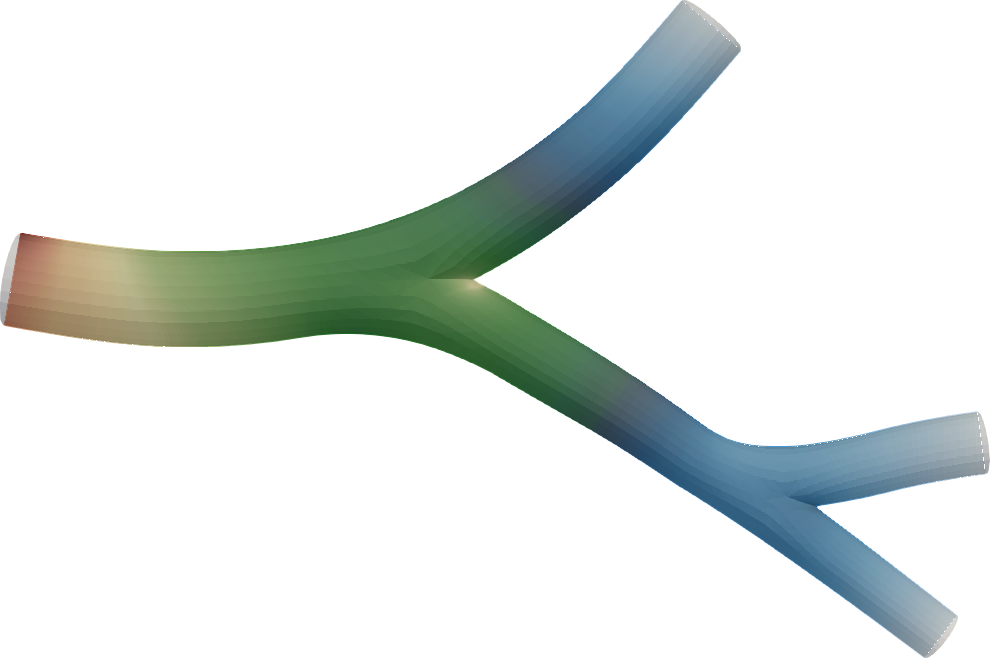}}
        \subfloat[Modes = $2$]{\label{fig:pmodes2}\includegraphics[width=.33\linewidth]{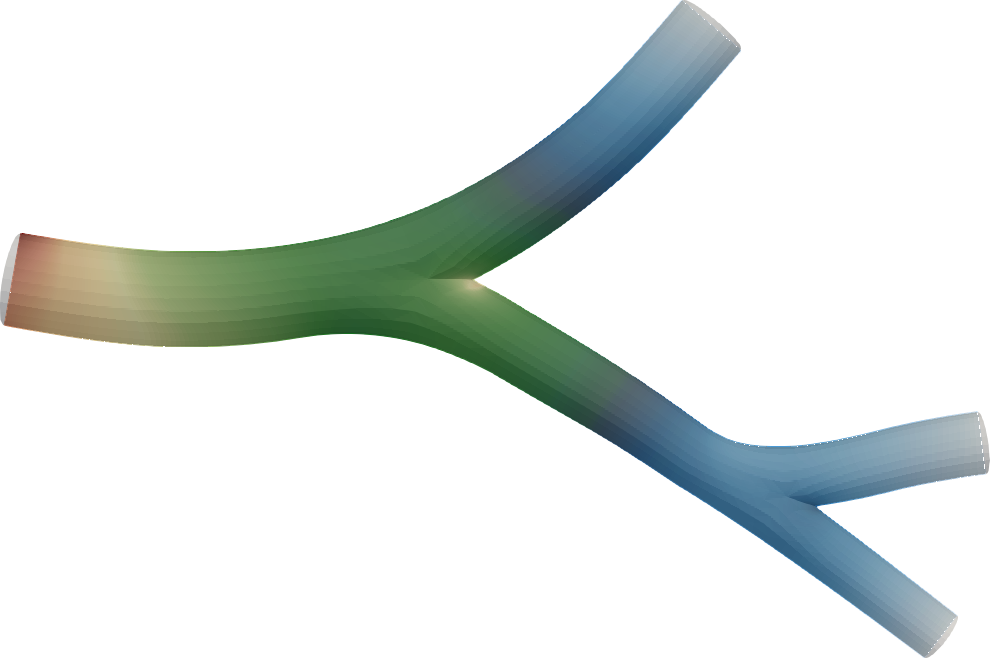}}
        \subfloat[Modes = $3$]{\label{fig:pmodes3}\includegraphics[width=.33\linewidth]{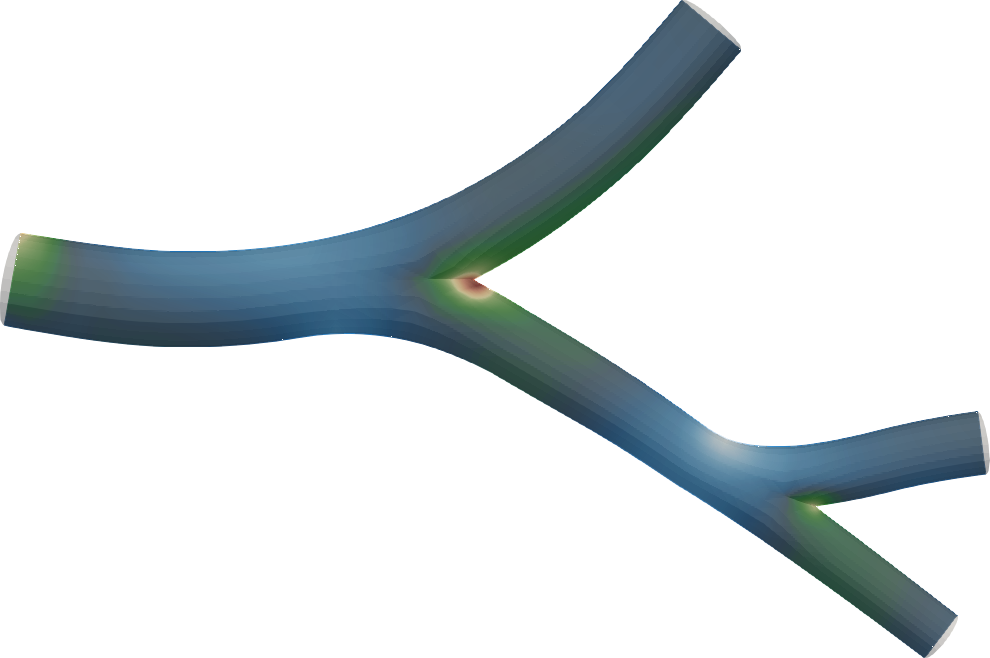}}
    \end{subfigure}

    \vspace{0.3cm} 

    \begin{subfigure}[b]{1\textwidth}
        \centering
        \subfloat[Modes = $4$]{\label{fig:pmodes4}\includegraphics[width=.33\linewidth]{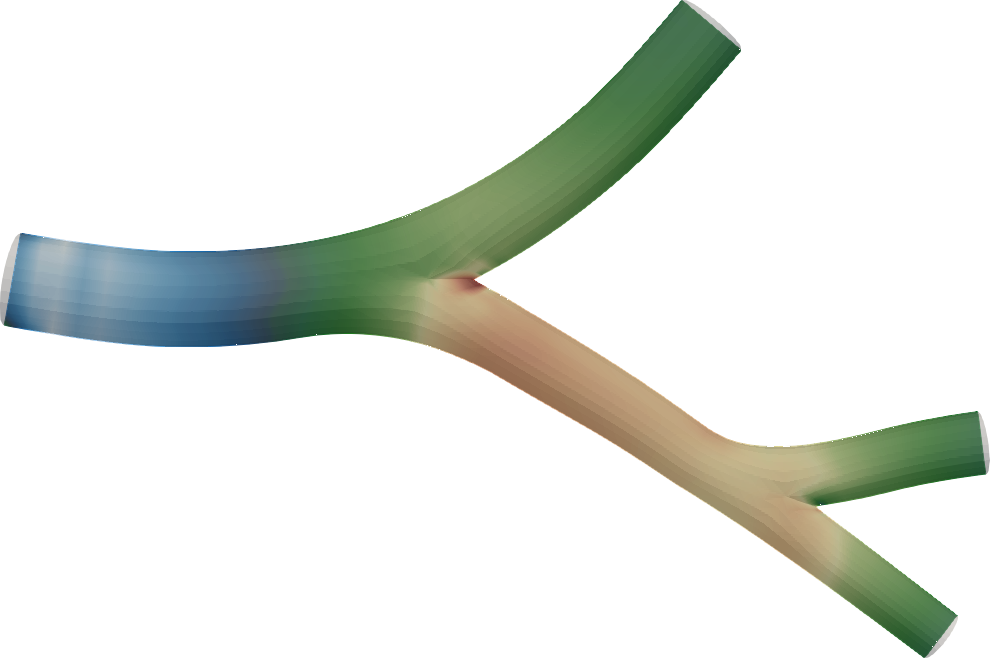}}
        \subfloat[Modes = $5$]{\label{fig:pmodes5}\includegraphics[width=.33\linewidth]{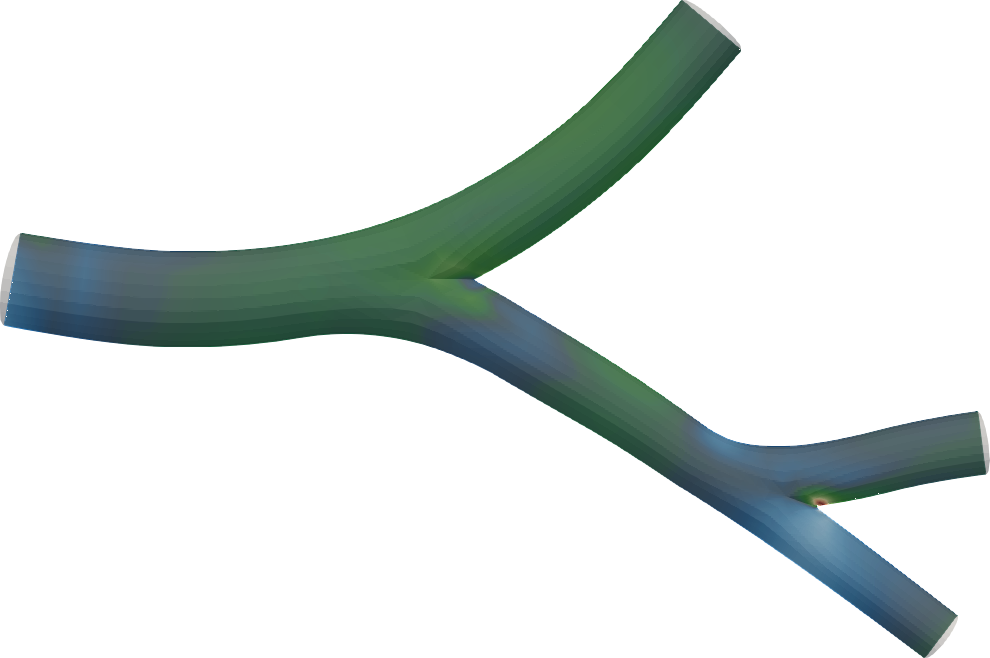}}
        \subfloat[Modes = $6$]{\label{fig:pmodes6}\includegraphics[width=.33\linewidth]{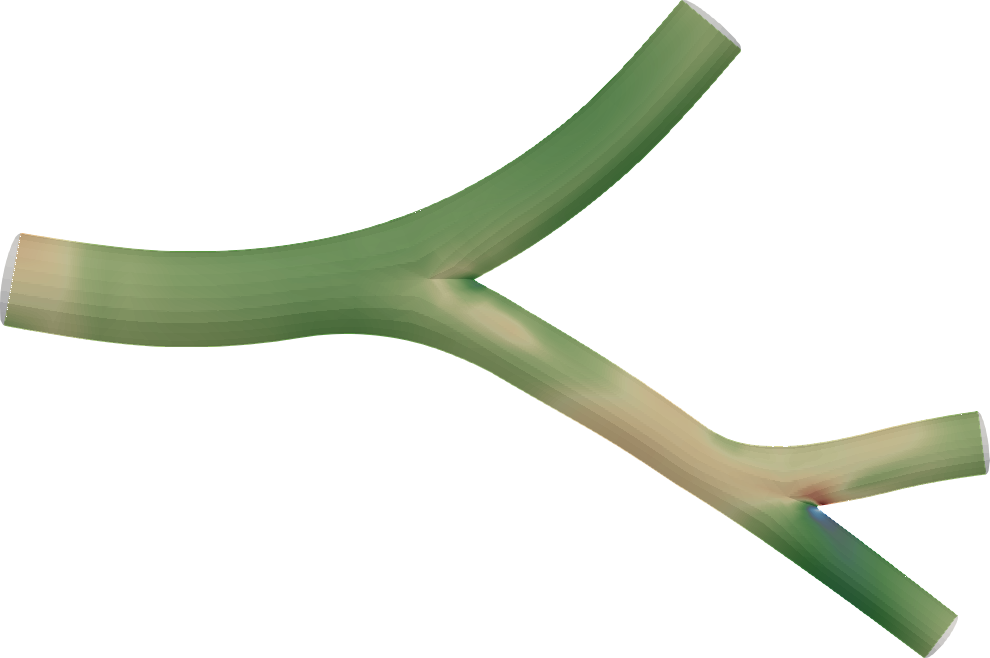}}
    \end{subfigure}

    \begin{subfigure}[b]{1\textwidth}
        \centering
        {\label{fig:pscale1}\includegraphics[width=0.25\linewidth]{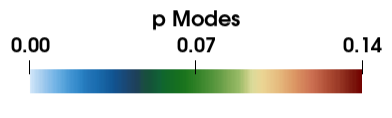}}          
    \end{subfigure}

    \medskip
    \caption{First six POD modes of the pressure field for the cerebrovascular bifurcation geometry. The modes are ordered according to decreasing energetic contribution, with the first mode representing the dominant pressure distribution.
}
    \label{fig:POD_P_modes}
\end{figure}

A POD decomposition is performed on the pressure field. The resulting pressure modes, illustrated in \fig{fig:POD_P_modes}, show a similar hierarchical structure, with the leading modes showing the dominant pressure distribution and global pressure gradients across the bifurcation, while higher-order modes indicate localized pressure variations. 

The obtained velocity and pressure POD modes are the foundation of the two reduced-order modeling  approaches examined in this study. For both models, the $3$ POD modes are retained for the velocity and pressure field to ensure a fair comparison between the POD-G and POD-RC approaches. In other words, both models rely on the same reduced representation of the full-order solution.

\subsection{POD-Galerkin ROM}
\label{sec:POD-G}

The POD–Galerkin projection-based ROM (POD-G) does not rely on an explicit split between training and test data. After projecting the Navier-Stokes equations onto the POD velocity and pressure modes, as described in Section~\ref{sec:intrusive_pod}, the reduced system is obtained from the full-order snapshots. The obtained reduced dynamical equations are integrated from the initial condition and then moved forward over the entire simulation time horizon. In our study, the test inlet velocity signal is explicitly defined from the beginning of the simulation with no segment of the test signal used for model calibration or initialization. Consequently, the POD–G ROM performs a fully physics-based online prediction for the entire duration of the test case (24 s).

To quantitatively assess the model accuracy, we compute the $L^2$ error between FOM and ROM solutions as follows:

\begin{equation}
E_{\theta}(t) = 100 \cdot \dfrac{||\theta_{FOM}(t) - \theta_{ROM}(t)||_{L^2(\Omega)}}{||{\theta_{FOM}}(t)||_{L^2(\Omega)}},
\label{eq:l2Error}
\end{equation}

where $\theta$ denotes the variable of interest. 

\begin{figure}[htb!]
    \centering \includegraphics[width=0.5\textwidth]{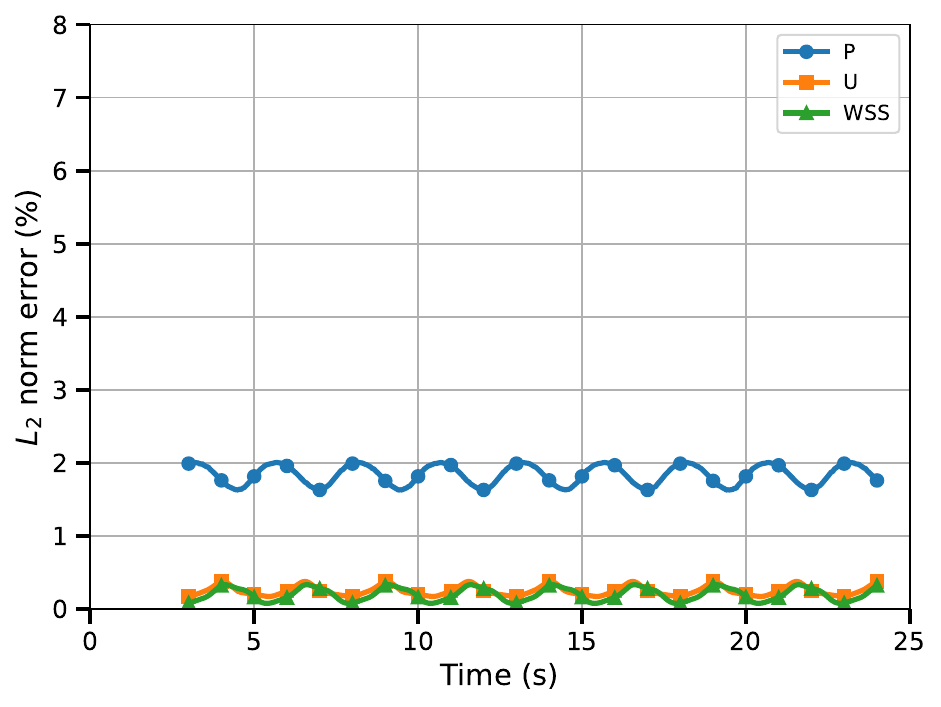}
        \caption{Time evolution  of the $L^2$ error \eqref{eq:l2Error} between FOM and POD-Galerkin projection ROM solutions for pressure (blue curve), velocity magnitude (orange curve) and wall shear stress (green curve). }   
        \label{fig:Rel_Error_online_All}
\end{figure}

\fig{fig:Rel_Error_online_All} shows the temporal evolution of the relative $L^22$ error for pressure ($p$), velocity magnitude $(U = |\mathbf{u}|)$ and wall shear stress $(\mathrm{WSS})$. The error evolution of all three variables follows the dynamics of the flow. They remain stable without any error growth and consistently small throughout the simulation time, confirming that POD-G ROM can efficiently and accurately predicts the flow dynamics.
Among the considered variables, it ${p}$ exhibits the highest error although it remains around $2\%$. ${U}$ shows significantly lower error, remaining below one percent for the entire simulation time and $\mathrm{WSS}$, although it is a derived variable from the velocity gradient, exhibits slightly lower error compared to the velocity field.

\begin{figure}[tb]
    \vspace{1cm}
    \centering
    \begin{overpic}[width=0.8\textwidth]{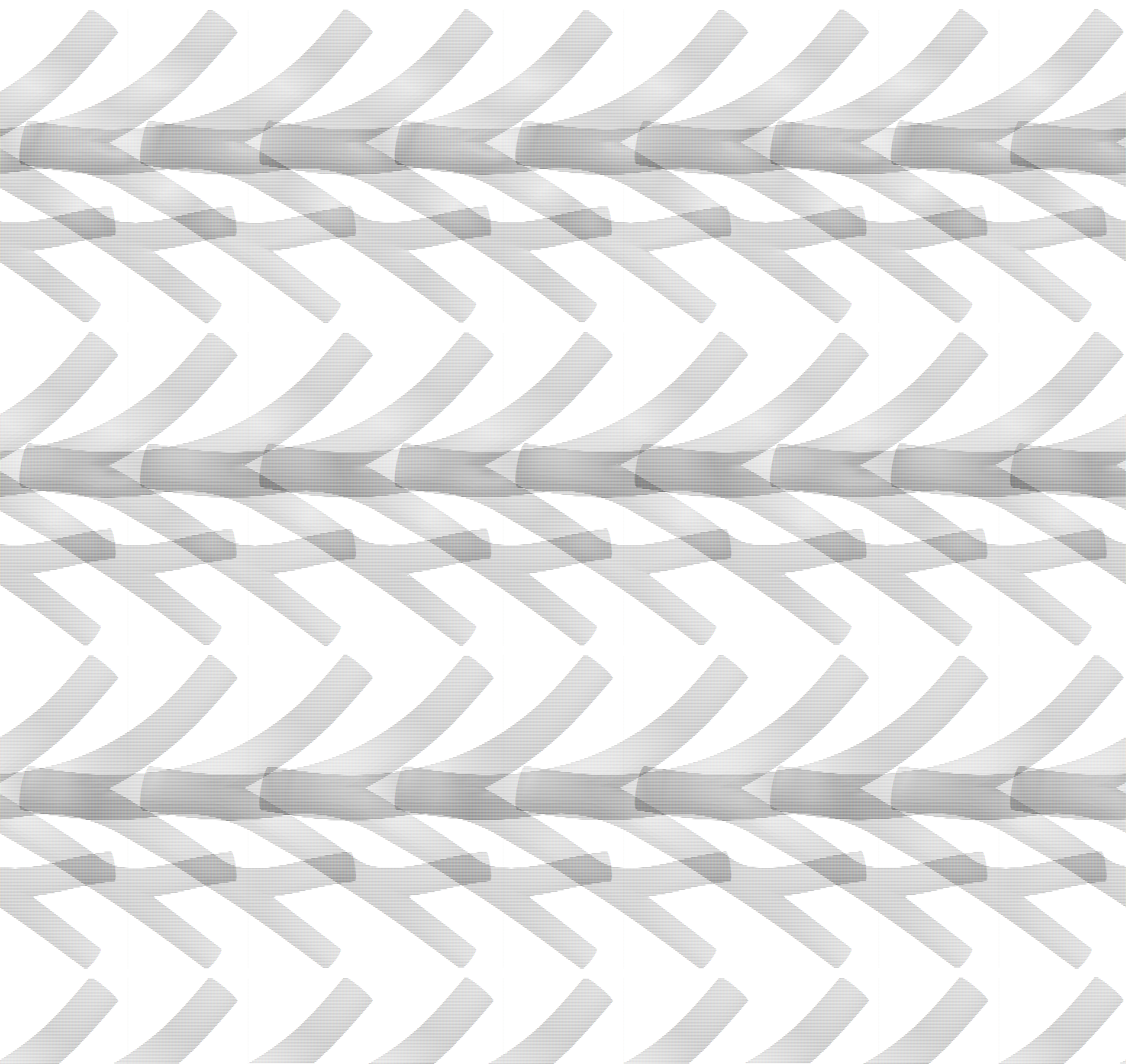}  
         \put(12,100.5){FOM}
         \put(45,100.5){POD-G}
         \put(77,100.5){$|\Delta p  |$}
         
         \put(10,92){\small{\textcolor{black}{$t$=4 s}}}
         \put(10,70){\small{\textcolor{black}{$t$=9 s}}}
         \put(10,49){\small{\textcolor{black}{$t$=14 s}}}
         \put(10,27){\small{\textcolor{black}{$t$=23 s}}}

         \put(42,92){\small{\textcolor{black}{$t$=4 s}}}
         \put(42,70){\small{\textcolor{black}{$t$=9 s}}}
         \put(42,49){\small{\textcolor{black}{$t$=14 s}}}
         \put(42,27){\small{\textcolor{black}{$t$=23 s}}}

         \put(74,92){\small{\textcolor{black}{$t$=4 s}}}
         \put(74,70){\small{\textcolor{black}{$t$=9 s}}}
         \put(74,49){\small{\textcolor{black}{$t$=14 s}}}
         \put(74,27){\small{\textcolor{black}{$t$=23 s}}}
        
      \end{overpic}
      \captionsetup{list=false}
  \caption{ Comparison of the evolution of $p$ given by the FOM (first column) and the POD–Galerkin ROM (second column) and their differences in absolute value (third column). }
   \label{fig:Vis_Online_P}
\end{figure}
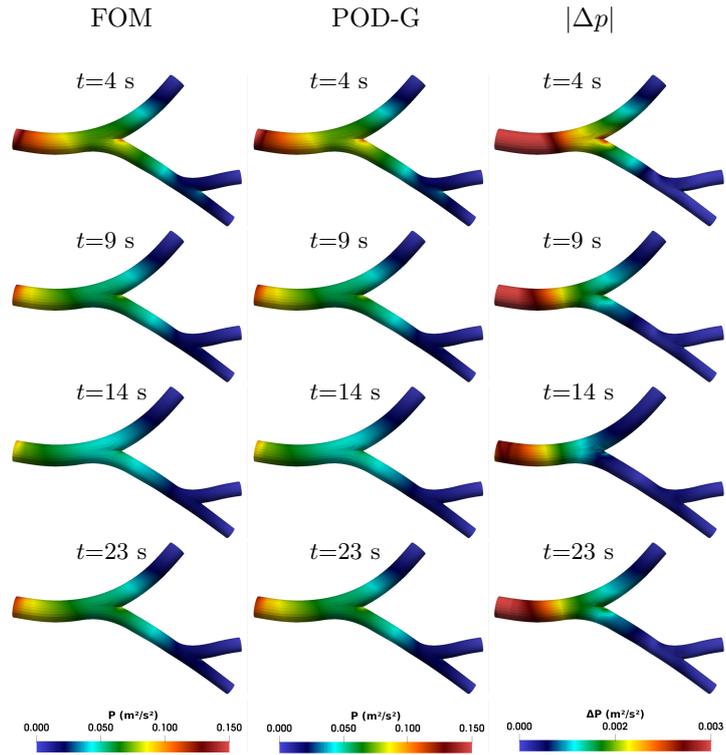

\fig{fig:Vis_Online_P} shows a qualitative comparison between FOM and POD-G ROM solutions for pressure field at four different time instances (i.e., $t = 4, 9, 14, 23 $ s) corresponding to various cycles of the test signal. As the test signal comprises multiple flow cycles, the reported snapshots are selected to represent the early, middle and last stages of the prediction horizon. This choice ensures that the observed agreement is not limited to a certain cycle. In general, ROM is capable of accurately reproducing the pressure field in the basilar artery bifurcation in term of spatial distribution and magnitude.   
The global pressure gradient from the inlet to the downstream branches is accurately represented at all reported times, along with the pressure redistribution that occurs in the bifurcation region. At all reported times, the difference in absolute pressure remains very small. The mismatch is observed mostly in the upstream of the bifurcation. The downstream branches, on the other hand, show almost no differences. Furthermore, the spatial distribution of $|\Delta p|$ remains uniform at all time instances demonstrating the absence of error accumulation throughout the prediction horizon. This trend is in agreement with the error evolution reported in \fig{fig:Rel_Error_online_All}.

\begin{figure}[htb!]
    \vspace{1cm}
    \centering
    \begin{overpic}[width=0.8\textwidth]{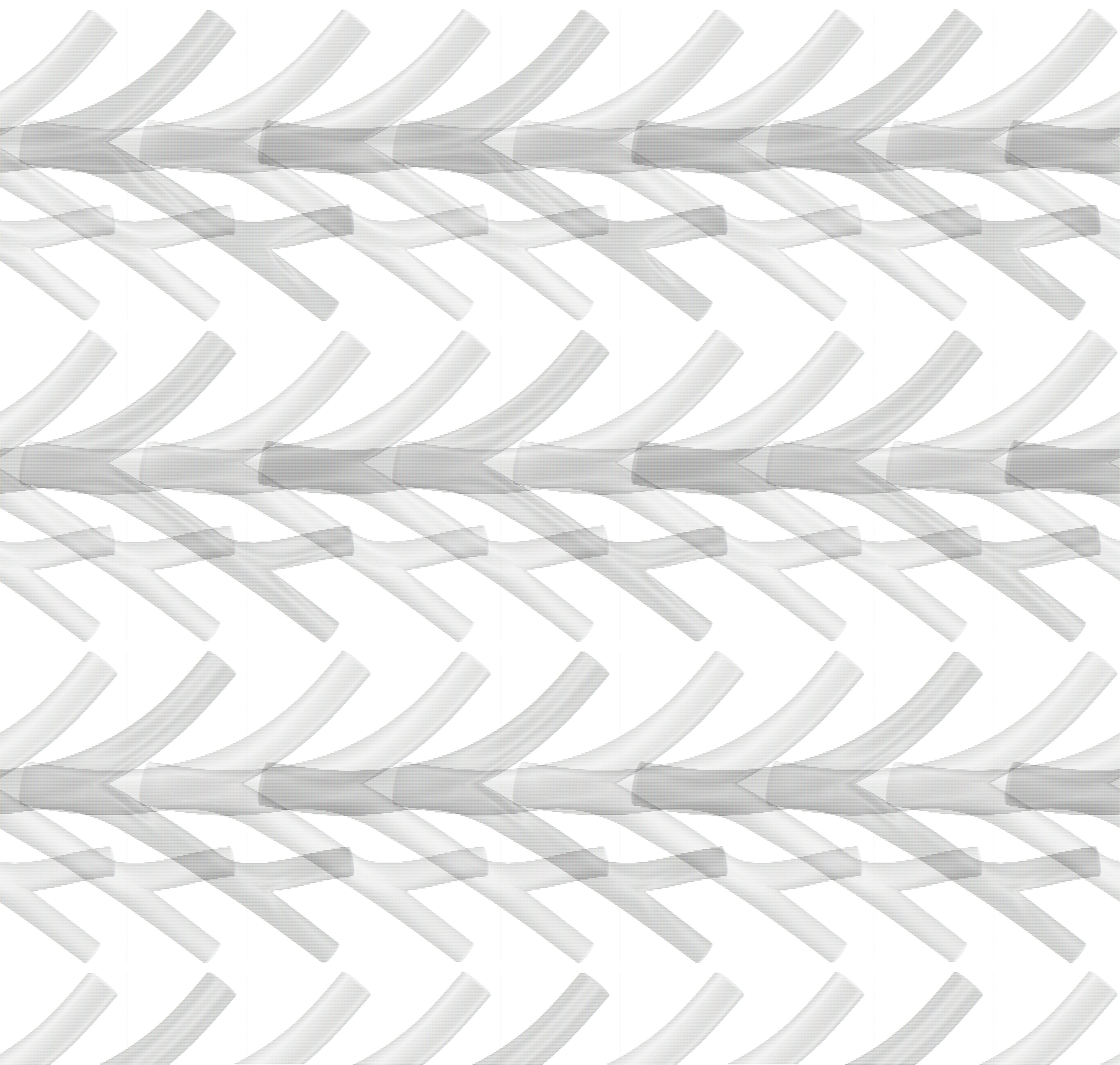}  
         \put(12,100.5){FOM}
         \put(45,100.5){POD-G}
         \put(77,100.5){$|\Delta U  |$}
         
         \put(10,92){\small{\textcolor{black}{$t$=4 s}}}
         \put(10,70){\small{\textcolor{black}{$t$=9 s}}}
         \put(10,49){\small{\textcolor{black}{$t$=14 s}}}
         \put(10,27){\small{\textcolor{black}{$t$=23 s}}}

         \put(42,92){\small{\textcolor{black}{$t$=4 s}}}
         \put(42,70){\small{\textcolor{black}{$t$=9 s}}}
         \put(42,49){\small{\textcolor{black}{$t$=14 s}}}
         \put(42,27){\small{\textcolor{black}{$t$=23 s}}}

         \put(74,92){\small{\textcolor{black}{$t$=4 s}}}
         \put(74,70){\small{\textcolor{black}{$t$=9 s}}}
         \put(74,49){\small{\textcolor{black}{$t$=14 s}}}
         \put(74,27){\small{\textcolor{black}{$t$=23 s}}}
        
      \end{overpic}
      \captionsetup{list=false}
  \caption{ Comparison of the evolution of $U$ given by the FOM (first column) and the POD–Galerkin ROM (second column) and their differences in absolute value (third column). }
   \label{fig:Vis_Online_U}
\end{figure}

A qualitative comparison of the velocity magnitude between FOM and ROM is shown in \fig{fig:Vis_Online_U} at the same four time instances. ROM accurately captures the flow dynamics including high velocity region in the parent vessel and the flow redistribution at the daughter branches as well as the flow acceleration at bifurcation zone. 
The absolute velocity difference remains very small all over the domain. The mismatch is mainly observed near the wall and vicinity of bifurcation region where the velocity gradients are stronger while a negligible difference is observed in the core flow region.  It is evident that the velocity difference $|\Delta U|$, both the magnitude and spatial distribution, remain consistent without any error accumulation in time verifying the low relative error distribution reported in \fig{fig:Rel_Error_online_All}.

\begin{figure}[htb!]
    \vspace{1cm}
    \centering
    \begin{overpic}[width=0.8\textwidth]{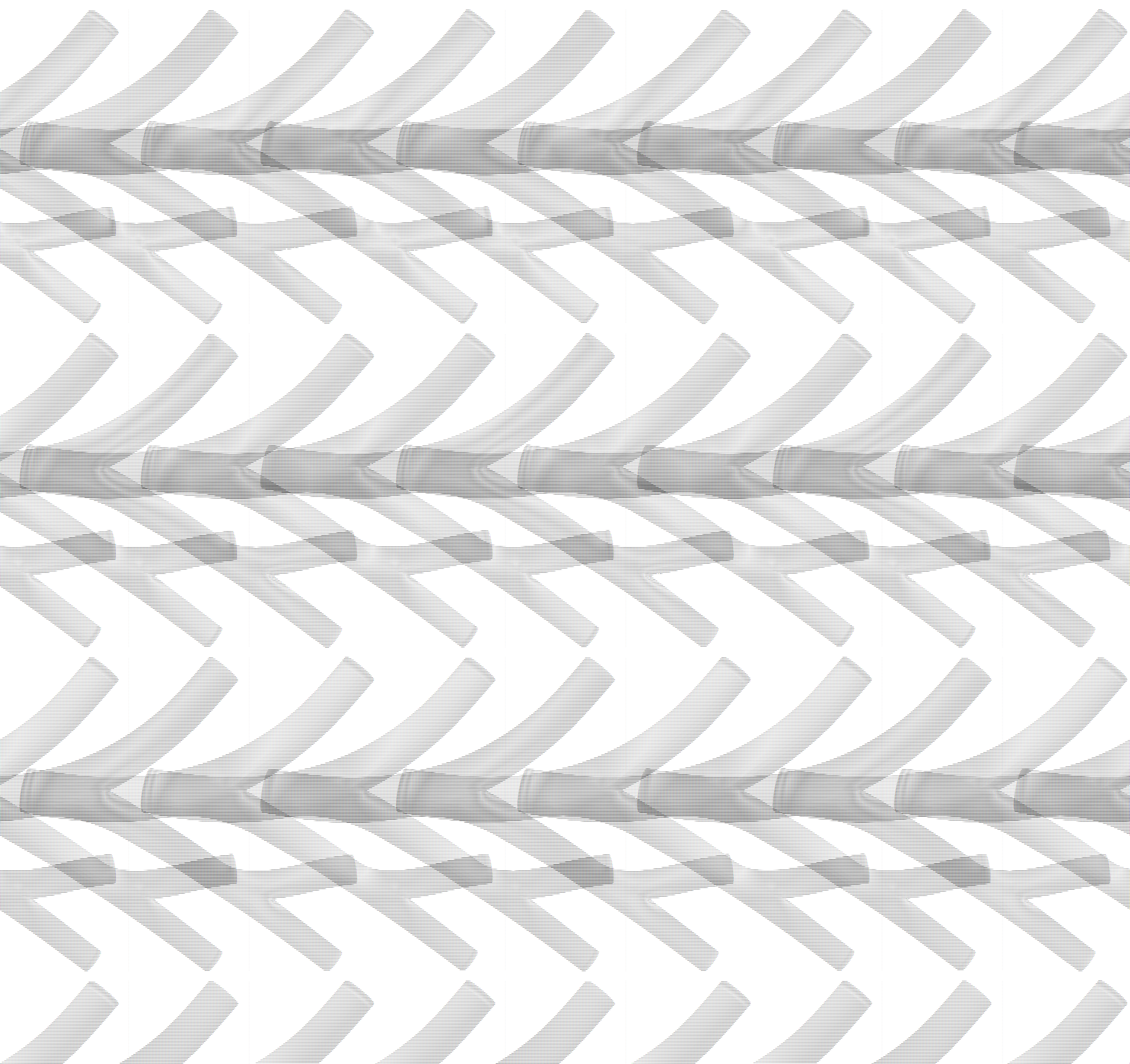}  
         \put(12,100.5){FOM}
         \put(45,100.5){POD-G}
         \put(74,100.5){$|\Delta\ WSS  |$}
         
         \put(10,92){\small{\textcolor{black}{$t$=4 s}}}
         \put(10,70){\small{\textcolor{black}{$t$=9 s}}}
         \put(10,49){\small{\textcolor{black}{$t$=14 s}}}
         \put(10,27){\small{\textcolor{black}{$t$=23 s}}}

         \put(42,92){\small{\textcolor{black}{$t$=4 s}}}
         \put(42,70){\small{\textcolor{black}{$t$=9 s}}}
         \put(42,49){\small{\textcolor{black}{$t$=14 s}}}
         \put(42,27){\small{\textcolor{black}{$t$=23 s}}}

         \put(74,92){\small{\textcolor{black}{$t$=4 s}}}
         \put(74,70){\small{\textcolor{black}{$t$=9 s}}}
         \put(74,49){\small{\textcolor{black}{$t$=14 s}}}
         \put(74,27){\small{\textcolor{black}{$t$=23 s}}}
        
      \end{overpic}
      \captionsetup{list=false}
  \caption{ Comparison of the evolution of $\mathrm{WSS}$
   given by the FOM (first column) and the POD–Galerkin ROM (second column) and their differences in absolute value (third column). }
   \label{fig:Vis_Online_WSS}
\end{figure}

Finally, the wall shear stress comparison in time and their differences are depicted in \fig{fig:Vis_Online_WSS} for the same time instances. Although $\mathrm{WSS}$ is a derived quantity depending on the velocity gradients, ROM is able to predict the main spatial pattern including the region of higher shear stress near the wall and bifurcation zone.
The absolute WSS differences are still small and only happen in certain places, mostly near bifurcation region where gradient-based quantities are most sensitive. The $\mathrm{WSS}$ error does not increase over time demonstrating  that the ROM prediction is stable and accurate.

It is worth mentioning that the absolute error field, shown in the right column of \fig{fig:Vis_Online_P}-\fig{fig:Vis_Online_WSS}, is reported using a different color bar scale than the corresponding FOM and ROM solutions. This color bar is one to two orders of magnitude smaller than those of the original variables.  This clear difference in scales highlights the very small values of reconstruction error and further confirms the high fidelity of POD–G ROM.




\subsection{POD–Reservoir Computing ROM }
\label{sec:POD-RC}

In contrast to the POD-Galerkin projection-based approach, the POD-Reservoir Computing (POD-RC) ROM uses a data-driven prediction method based on the temporal evolution of POD coefficients. As explained in Section~\ref{sec:POD-RC}, the high-dimensional flow fields are initially projected onto a reduced POD space and the resulting modal coefficients are used to train the reservoir computing network.

The multi-harmonic inlet velocity signal, the blue curve shown in \fig{fig:training_test_signals}, is used during the training phase of the POD-RC ROM to learn the temporal evolution of the POD coefficients. A different inlet velocity signal is considered for the test phase, the orange curve in \fig{fig:training_test_signals}. As RC network is a type of recurrent neural network that requires an initial time context, a short part of the test signal $(3 s)$ is used to initialize the reservoir state. The trained POD-RC model then predicts the evolution of the POD coefficients over the remaining simulation time, corresponding to the subsequent 21 s.
The performance of the POD-RC ROM is evaluated through the following quantitative and qualitative comparisons of the reconstructed flow fields.

\begin{figure}[htb!]
    \centering \includegraphics[width=0.5\textwidth]{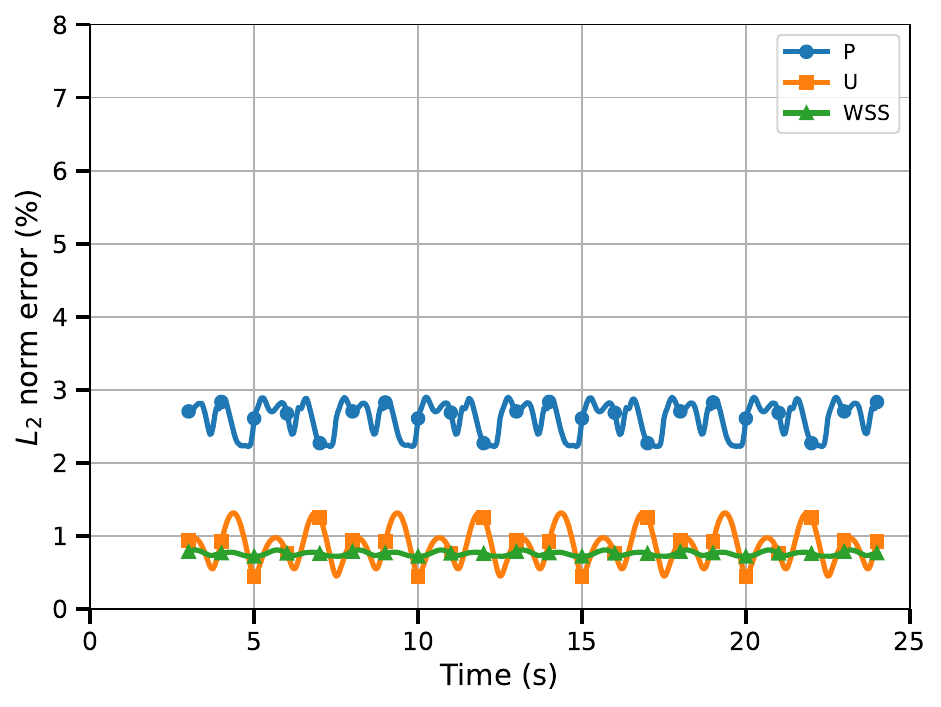}
        \caption{Time evolution  of the $L^2$ error \eqref{eq:l2Error} between FOM and POD-RC ROM solutions for pressure (blue curve), velocity magnitude (orange curve) and wall shear stress (green curve). }   
        \label{fig:Rel_Error_RC_All}
\end{figure}

We report on \fig{fig:Rel_Error_RC_All} the temporal evolution of the relative  $L^2$ error \ref{eq:l2Error} for the variable of interest in this study, namely p, U and $\mathrm{WSS}$. For all the variables,  the error remains low and stable without any error growth in time, confirming the robustness of POD-RC approach for this problem. Among the considered variables, the highest relative error is associated with the pressure which fluctuates approximately 3$\%$ over the entire simulation time. U and $\mathrm{WSS}$ exhibit significantly lower error, remaining below approximately 1.2$\%$.

\begin{figure}[htb!]
    \vspace{1cm}
    \centering
    \begin{overpic}[width=0.8\textwidth]{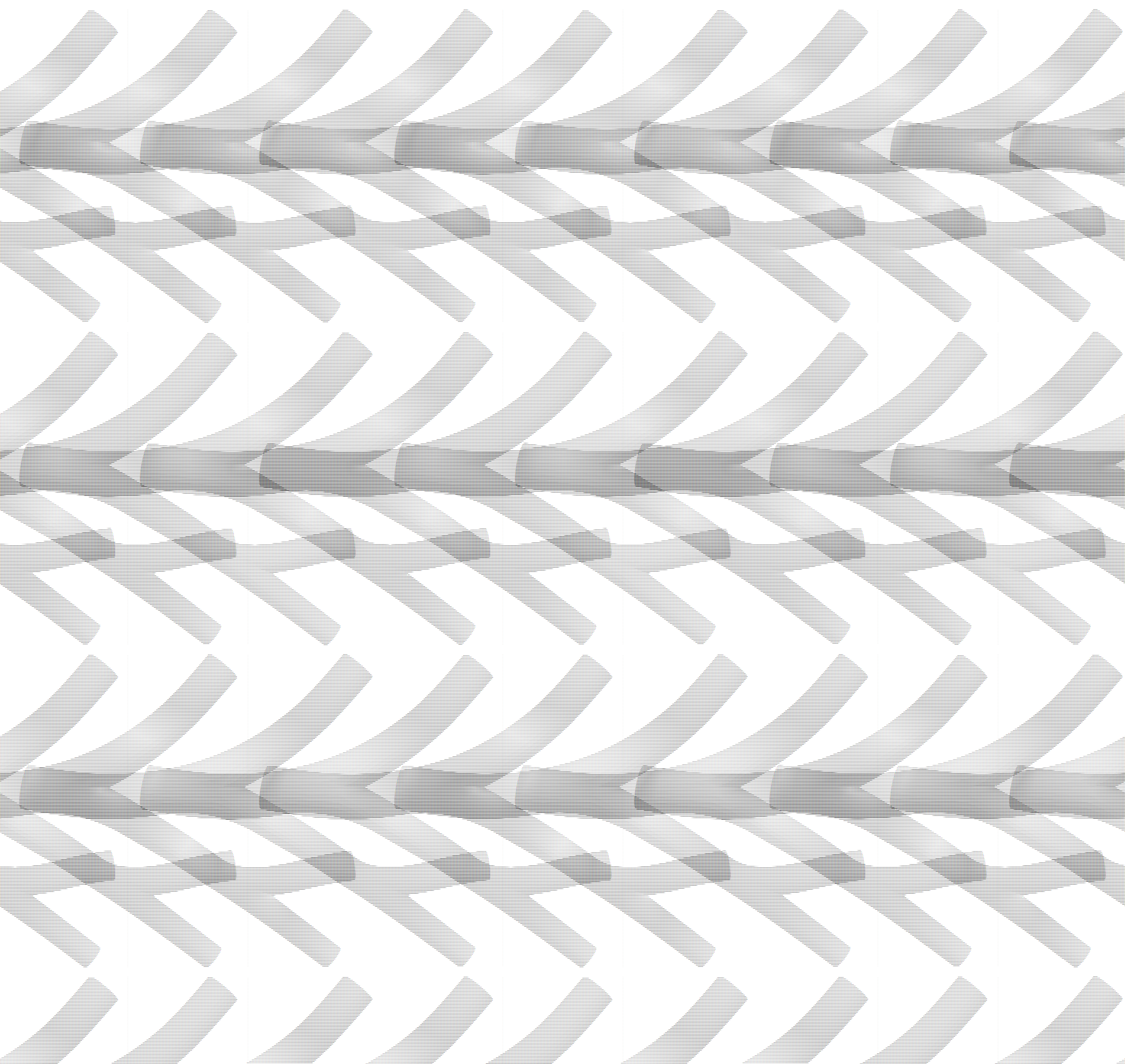}  
         \put(12,100.5){FOM}
         \put(45,100.5){POD-RC}
         \put(77,100.5){$|\Delta p  |$}
         
         \put(10,92){\small{\textcolor{black}{$t$=4 s}}}
         \put(10,70){\small{\textcolor{black}{$t$=9 s}}}
         \put(10,49){\small{\textcolor{black}{$t$=14 s}}}
         \put(10,27){\small{\textcolor{black}{$t$=23 s}}}

         \put(42,92){\small{\textcolor{black}{$t$=4 s}}}
         \put(42,70){\small{\textcolor{black}{$t$=9 s}}}
         \put(42,49){\small{\textcolor{black}{$t$=14 s}}}
         \put(42,27){\small{\textcolor{black}{$t$=23 s}}}

         \put(74,92){\small{\textcolor{black}{$t$=4 s}}}
         \put(74,70){\small{\textcolor{black}{$t$=9 s}}}
         \put(74,49){\small{\textcolor{black}{$t$=14 s}}}
         \put(74,27){\small{\textcolor{black}{$t$=23 s}}}
        
      \end{overpic}
      \captionsetup{list=false}
  \caption{ Comparison of the evolution of p
   given by the FOM (first column) and the POD–Reservoir Computing ROM (second column) and their differences in absolute value (third column). }
   \label{fig:Vis_RC_P}
\end{figure}

\begin{figure}[htb!]
    \vspace{1cm}
    \centering
    \begin{overpic}[width=0.8\textwidth]{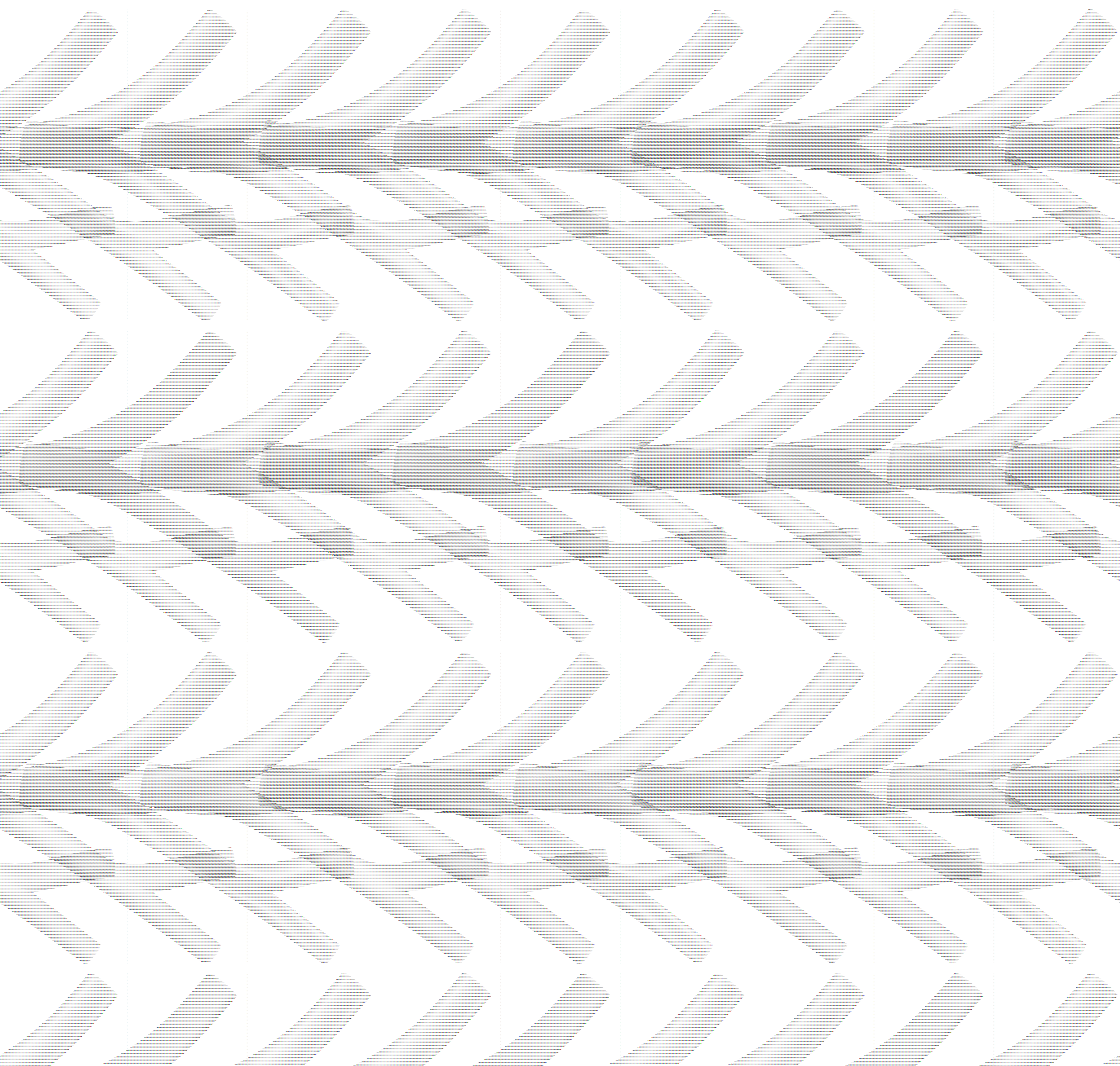}  
         \put(12,100.5){FOM}
         \put(45,100.5){POD-RC}
         \put(77,100.5){$|\Delta U  |$}
         
         \put(10,92){\small{\textcolor{black}{$t$=4 s}}}
         \put(10,70){\small{\textcolor{black}{$t$=9 s}}}
         \put(10,49){\small{\textcolor{black}{$t$=14 s}}}
         \put(10,27){\small{\textcolor{black}{$t$=23 s}}}

         \put(42,92){\small{\textcolor{black}{$t$=4 s}}}
         \put(42,70){\small{\textcolor{black}{$t$=9 s}}}
         \put(42,49){\small{\textcolor{black}{$t$=14 s}}}
         \put(42,27){\small{\textcolor{black}{$t$=23 s}}}

         \put(74,92){\small{\textcolor{black}{$t$=4 s}}}
         \put(74,70){\small{\textcolor{black}{$t$=9 s}}}
         \put(74,49){\small{\textcolor{black}{$t$=14 s}}}
         \put(74,27){\small{\textcolor{black}{$t$=23 s}}}
        
      \end{overpic}
      \captionsetup{list=false}
  \caption{ Comparison of the evolution of U
   given by the FOM (first column) and the POD–Reservoir Computing ROM (second column) and their differences in absolute value (third column). }
   \label{fig:Vis_RC_U}
\end{figure}

\begin{figure}[htb!]
    \vspace{1cm}
    \centering
    \begin{overpic}[width=0.9\textwidth]{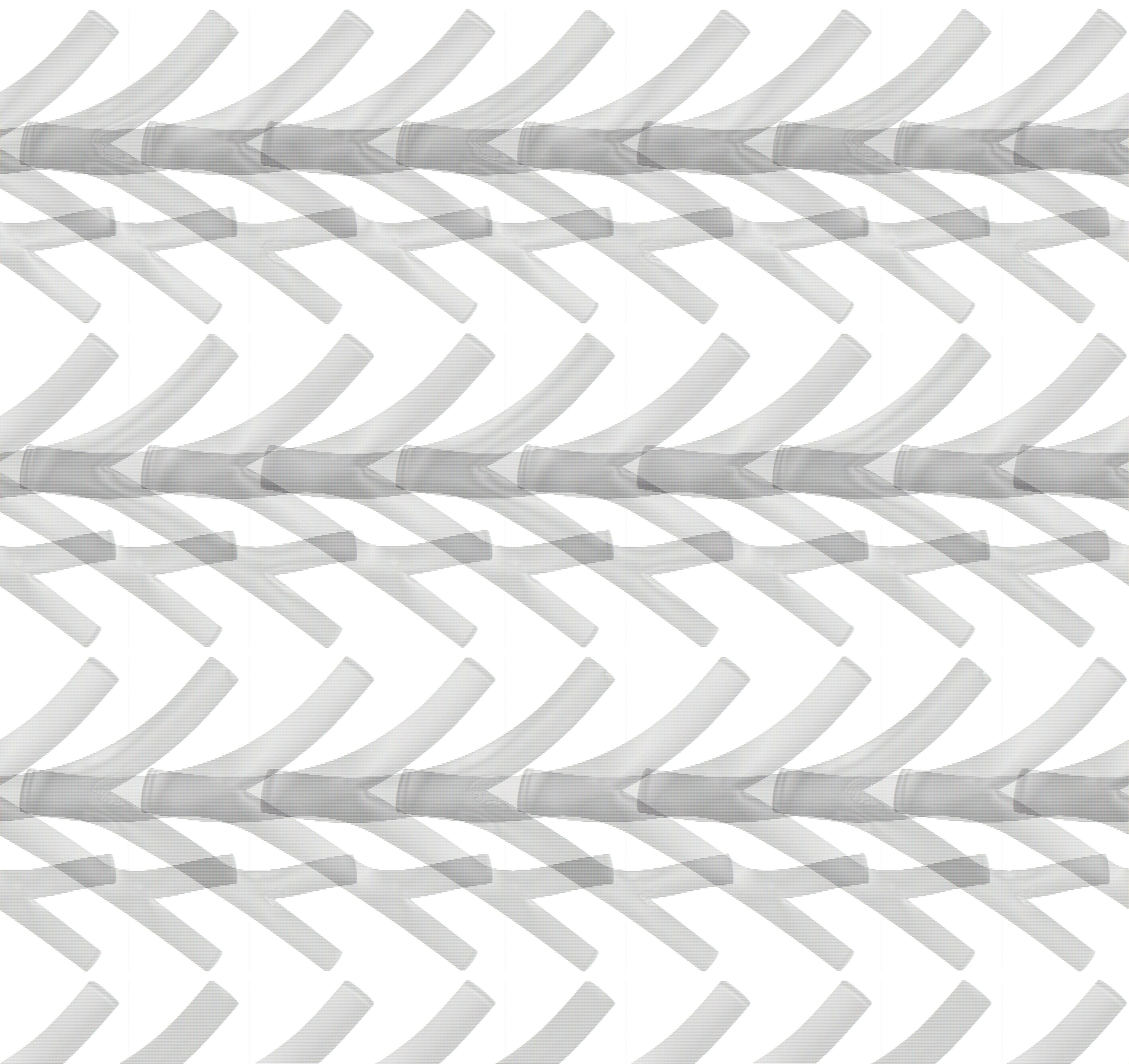}  
         \put(12,100.5){FOM}
         \put(45,100.5){POD-RC}
         \put(74,100.5){$|\Delta \ WSS  |$}
         
         \put(10,92){\small{\textcolor{black}{$t$=4 s}}}
         \put(10,70){\small{\textcolor{black}{$t$=9 s}}}
         \put(10,49){\small{\textcolor{black}{$t$=14 s}}}
         \put(10,27){\small{\textcolor{black}{$t$=23 s}}}

         \put(42,92){\small{\textcolor{black}{$t$=4 s}}}
         \put(42,70){\small{\textcolor{black}{$t$=9 s}}}
         \put(42,49){\small{\textcolor{black}{$t$=14 s}}}
         \put(42,27){\small{\textcolor{black}{$t$=23 s}}}

         \put(74,92){\small{\textcolor{black}{$t$=4 s}}}
         \put(74,70){\small{\textcolor{black}{$t$=9 s}}}
         \put(74,49){\small{\textcolor{black}{$t$=14 s}}}
         \put(74,27){\small{\textcolor{black}{$t$=23 s}}}
        
      \end{overpic}
      \captionsetup{list=false}
  \caption{ Comparison of the evolution of $WSS$
   given by the FOM (first column) and the POD–Reservoir Computing ROM (second column) and their differences in absolute value (third column). }
   \label{fig:Vis_RC_WSS}
\end{figure}

Figure~\ref{fig:Vis_RC_P} illustrates the qualitative comparison between the FOM and POD–RC pressure fields at four different time intervals: $t = 4, 9, 14, 23$~s, which correspond to various oscillation cycles of the test signal. The POD–RC model can accurately reconstruct the overall pressure distribution in the basilar artery bifurcation, capturing both the spatial structure and the temporal evolution of the pressure field over the prediction horizon. 
The absolute pressure difference remains small over the entire domain. The differences are mostly observed at the upstream part of the parent vessel, while the downstream branches show negligible differences. Importantly, the spatial distribution and magnitude of $|\Delta p|$ remain consistent across all selected time instances demonstrating that the POD-RC prediction does not exhibit error accumulation over time. This behavior aligns with the bounded pressure $L^2$ error reported in the quantitative analysis.

Figure~\ref{fig:Vis_RC_U} shows the qualitative comparison of the velocity magnitude between the FOM and POD-RC predictions at the same four time instances. The POD–RC method accurately reproduces the main flow features including the high-velocity area in the parent vessel and the redistribution of the flow into the daughter branches. The agreement between the FOM and POD-RC solutions stays the same for the entire prediction interval.
The absolute velocity difference $|\Delta U|$ stays small and is mostly localized near the walls of the vessel, where the velocity gradients are strongest. The magnitude of the velocity error remains consistent across all selected time instances confirming the velocity reconstruction does not deteriorate over time.

Figure~\ref{fig:Vis_RC_WSS} presents a qualitative comparison of the wall shear stress from the FOM and the POD-RC model, along with the absolute differences.  The POD–RC method does a good job of reproducing the main spatial patterns of wall shear stress, such as areas of high values near the bifurcation and along the vessel walls.
The absolute $\mathrm{WSS}$ differences are still small over most of the wall surface. They are mostly localized in small areas near the bifurcation where gradient-based quantities are naturally more sensitive to small changes in velocity. As observed for pressure and velocity, the spatial distribution of $|\Delta \mathrm{WSS}|$ remains small.

\subsection{Comparative Assessment of Reduced Order Models}
\label{sec:comparison}

\begin{figure}[t]
    \centering

    \begin{subfigure}{0.48\textwidth}
        \centering
        \includegraphics[width=\linewidth]{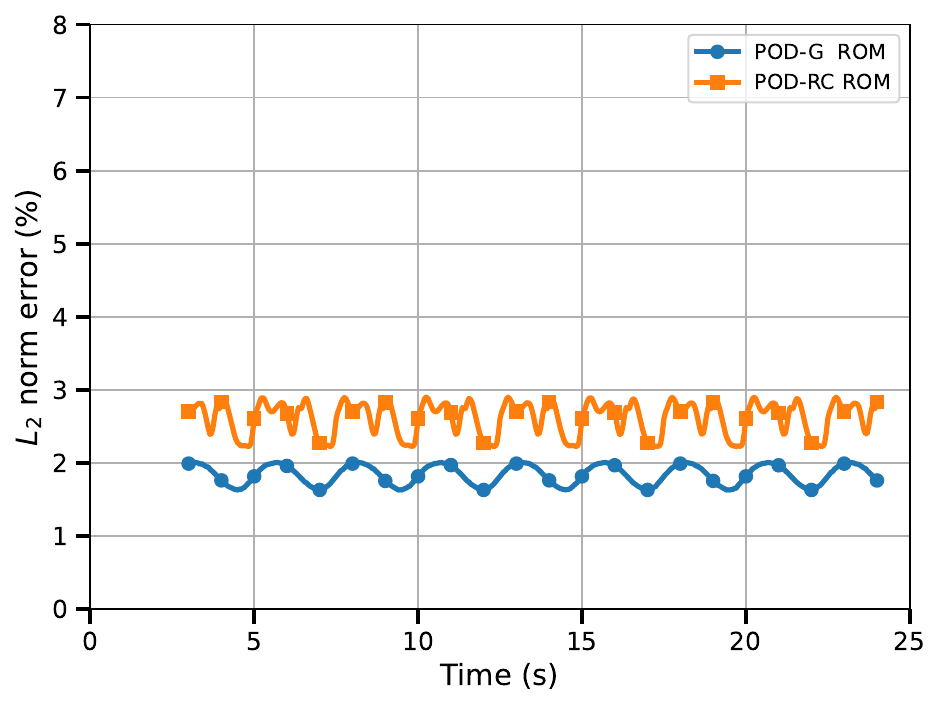}
        \caption{Pressure (p)}
        \label{subfig:Rel_Error_P}
    \end{subfigure}
    \hfill
    \begin{subfigure}{0.48\textwidth}
        \centering
        \includegraphics[width=\linewidth]{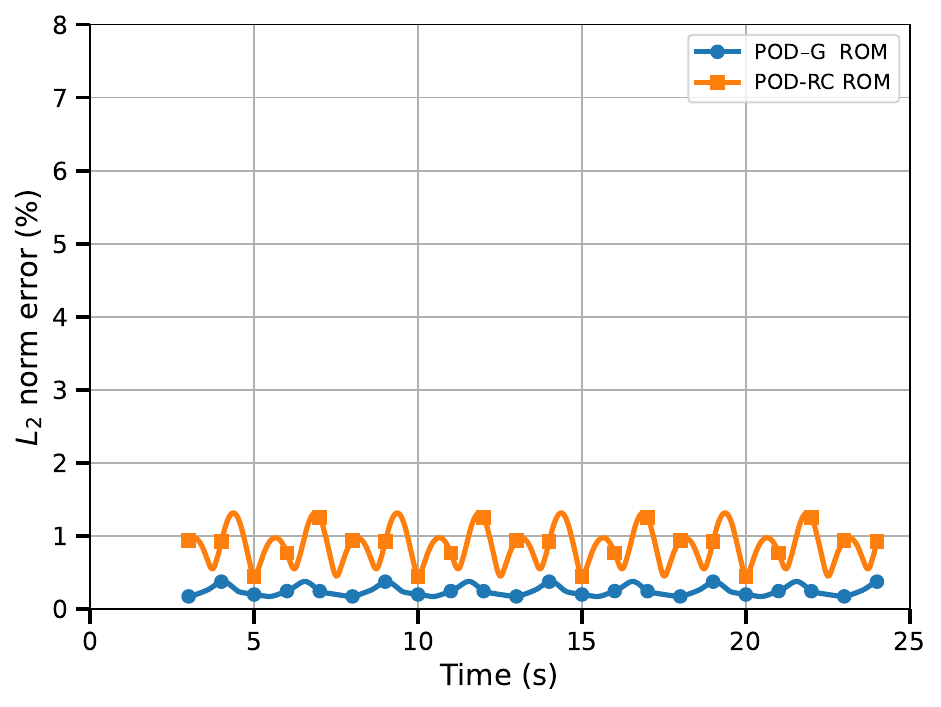}
        \caption{Velocity magnitude (U)}
        \label{subfig:Rel_Error_U}
    \end{subfigure}

    \vspace{0.3cm}

    \begin{subfigure}{0.48 \textwidth}
        \centering
        \includegraphics[width=\linewidth]{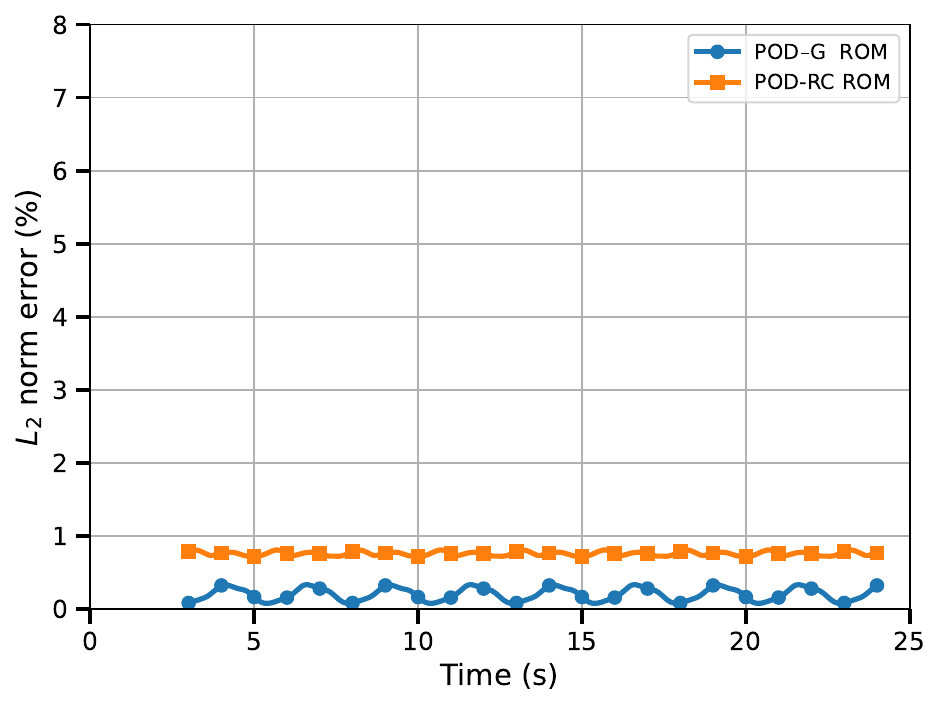}
        \caption{Wall Shear Stress ($\mathrm{WSS}$)}
        \label{subfig:Rel_Error_WSS}
    \end{subfigure}

    \caption{Relative L2 error comparison over time}
    \label{fig:Rel_Error_Comparison}
\end{figure}

In this section, a comparative analysis between the two proposed ROM approaches is presented: POD-G ROM (Intrusive ROM) and POD-RC ROM (Data-Driven); see \fig{fig:Rel_Error_Comparison}. This comparison is based on the temporal evolution of $L^2$ error for pressure \fig{subfig:Rel_Error_P}, velocity magnitude \fig{subfig:Rel_Error_U}, and wall shear stress \fig{subfig:Rel_Error_WSS} to provide insights into the accuracy and robustness of these two models.


For all three quantities, both models exhibit very small and stable error levels throughout the prediction horizon, as also observed in the previous section. This behavior confirms that there are no instabilities in the POD-reduced system, and both models are capable of accurately and efficiently predicting the flow dynamics over time.

For the pressure field (\fig{subfig:Rel_Error_P}), the POD-G ROM consistently achieves lower error levels than the POD-RC ROM. It fluctuates around 1.6$\%$ to 2.0$\%$ while the data-driven approach oscillates between higher error values, $2.3\%$--$3.0\%$ 
A similar trend is observed for velocity magnitude (\fig{subfig:Rel_Error_U}); The POD-G ROM exhibits very low error, below 0.4$\%$, while the POD-RC ROM reaches values up to 1.2$\%$.
In the case of wall shear stress (\fig{subfig:Rel_Error_WSS}), both models produce low relative errors. However, the POD-RC ROM again shows higher values compared to the lower errors obtained from the projection-based method.

In conclusion, POD-G ROM outperforms the POD-RC ROM model in predicting all considered quantities as it explicitly incorporates the governing physical equations into the projection-based formulations. 
On the other hand, the POD-RC model relies only on the data-driven temporal modeling of the reduced coefficients, leading to slightly higher errors but still gives stable and accurate long-term predictions.

\subsection{Computational efficiency}
\label{sec:computational_cost}

All simulations were performed on a workstation with an 11th Gen Intel(R) Core(TM) i7-11700 CPU running at 2.50 GHz and 32 GB of RAM. A single FOM simulation on this platform takes about XX seconds to complete, highlighting the high computational cost of repeating high-fidelity simulations.

In reduced-order modeling, the computational cost is divided into two parts: offline and online. In the POD-G ROM, the offline phase consists of calculating the POD modes and using intrusive projection to reconstruct the reduced dynamical system. This offline phase requires a significant amount of computer cost and takes approximately XX seconds. The POD-RC ROM, on the other hand, has a significantly lower offline cost. Once the POD coefficients have been obtained, the reservoir computing network can be trained quickly without requiring intrusive access to the governing equations. The POD-RC method requires approximately $12000$ seconds of offline training time, which is significantly lower than that required by POD-G ROM method.

After the offline phase is over, both reduced-order models have a very fast online phase. The POD–G ROM advances the reduced dynamical system over the entire prediction horizon in about $0.1$ seconds, while the POD–RC ROM predicts the future POD coefficients and reconstructs the corresponding  physical fields in about $10$ seconds. In both cases, the required time for online prediction is several orders of magnitude smaller than that required for a FOM simulation. Overall, both reduced-order models make it possible to make accurate online predictions; however, the POD-RC ROM is more efficient due to its significantly shorter offline cost.

\begin{table}[htbp]
\centering
\begin{tabular}{lccc}
\hline
\textbf{Method} & \textbf{Training Time (s)} & \textbf{Online Time (s)} & \textbf{Speedup vs FOM} \\ \hline
High-Fidelity CFD & $\sim$ 12,000 & $\sim$ 12,000 & 1.0 \\
POD--Galerkin & $\sim$ 15 & $\sim$ 10 & $\sim$ 1,000 \\
Reservoir Computing & $\sim$ 2 & $\sim$ 0.1 & $\sim$ 100,000\\ \hline
\end{tabular}

\caption{Comparison of computational costs. The data-driven surrogates (RNN and RC) offer orders-of-magnitude speedups compared to the full-order solver.}
\label{tab:cpu_time}
\end{table}

\section{Concluding Remarks}

In this study, we examined and compared two reduced-order models for unsteady cerebrovascular hemodynamics: a physics-based POD-Galerkin projection ROM and a data-driven POD-Reservoir Computing (POD-RC) ROM. Both methods were constructed from high-fidelity CFD simulations of an idealized bifurcation of the basilar artery. They were examined in terms of accuracy, computational efficiency, and temporal stability in predicting: pressure (p), velocity magnitude (U), and wall shear stress ($\mathrm{WSS}$).

The projection-based POD-G ROM shows a promising approach, providing accurate and stable predictions throughout the entire simulation time. It consistently produced the lowest relative error for all the considered quantities, benefiting from the explicit incorporation of the governing Navier–Stokes equations into the reduced system. The qualitative comparison further confirms its capability in reconstructing the spatial patterns and temporal evolution of p, U and $\mathrm{WSS}$ over time for all the flow cycles without any error growth. The POD-RC ROM, although it is a pure data-driven approach, is capable of accurately and consistently predicting long-time frames. This approach successfully captured the main flow features and maintained the error levels limited throughout the simulation, although it exhibited slightly higher error levels than the POD–G ROM. Notably, its performance was highlighted, as it can learn complex reduced-order dynamics with minimal training and without direct access to the governing equations.

From a computational point of view, both reduced-order models achieved  notable speedup compared to  the full-order CFD simulations. However, the POD–RC ROM is significantly faster in the offline phase, as it avoids the expensive intrusive projection step that is required by the POD–G method. These characteristics make POD-RC an attractive approach for problems in which access to the governing equations is limited, the underlying dynamic is very complex or the main object of ROM is to construct a model for fast and flexible prediction.

In general, the results show that the POD-G ROM is the most accurate choice as it is a physics-based approach. The POD-RC ROM, on the other hand, is a good balance between accuracy and efficiency. The latter is especially promising for cerebrovascular applications that need to work in real time or be specific to each patient where the fast surrogate modeling with reliable long-term predictions is essential. Future work will focus on extending these methods to more complex geometries, the introduction of non-Newtonian rheology, turbulence models, more appropriate pressure outlets in the governing physics and a wider range of parameters for the consideration of the reduced order model. The present work lays the foundation for combining Reservoir Computation with governing physics in full-order and reduced order form and offers promise for faster physics-informed neural networks.

\section{Acknowledgments}
We acknowledge the support provided by the European Research Council Executive Agency by the Consolidator Grant project AROMA-CFD "Advanced Reduced Order Methods with Applications in Computational Fluid Dynamics" - GA 681447, H2020-ERC CoG 2015 AROMA-CFD, the European Union’s Horizon 2020 research and innovation program under the Marie Skłodowska-Curie Actions, grant agreement 872442 (ARIA), PON “Research and Innovation on Green related issues” FSE REACT-EU 2021 project, PRIN NA FROM-PDEs project, INdAM-GNCS 2019-2021 projects.
\clearpage
\bibliographystyle{unsrt}
\bibliography{mybib}

\end{document}